\documentclass{article}

\usepackage[cp1250]{inputenc}
\usepackage[IL2]{fontenc}
\usepackage{yfonts,fancyhdr}
\usepackage{a4wide}
\usepackage{euscript}
\usepackage{amstext,amsbsy,amscd,amssymb}
\usepackage{amsmath,enumerate}
\usepackage{amsfonts}
\usepackage{mathrsfs,tensor,xy}
\usepackage{graphics}
\include{diagxy}

\long\def\proof #1{\noindent \emph{Proof.}\ #1 \hfill $\spadesuit$
\medskip}

\newcounter{numdef}
\long\def\definition #1 {\refstepcounter{numdef} \noindent {\bf
Definition \arabic{numdef}.} #1

\medskip}

\newcounter{numthe}
\long\def\theorem #1{\refstepcounter{numthe} \noindent {\bf Theorem
\arabic{numthe}.} #1

\medskip}

\newcounter{numlem}
\long\def\lemma #1{\refstepcounter{numlem}  \noindent {\bf Lemma
\arabic{numlem}.} #1

\medskip}

\long\def\example{\noindent {\bf Example.}\ }

\long\def\remark{\noindent {\bf Remark.}\ }

\def\Diff{\mathop {\rm Diff} \nolimits}

\def\End{\mathop {\rm End} \nolimits}

\def\Aut{\mathop {\rm Aut} \nolimits}

\def\ind{\mathop {\rm ind} \nolimits}

\def\ad{\mathop {\rm ad} \nolimits}
\def\Ad{\mathop {\rm Ad} \nolimits}

\def\esssupp{\mathop {\rm ess\,supp}\limits}
\def\sign{\mathop {\rm sign} \nolimits}
\def\Gau{\mathop {\rm Gau} \nolimits}
\def\gau{\mathop \mfrak{gau} \nolimits}
\def\im{\mathop {\rm im} \nolimits}
\def\tr{\mathop {\rm tr} \nolimits}
\def\rk{\mathop {\rm rk} \nolimits}
\def\vol{\mathop {\rm vol} \nolimits}
\def\Der{\mathop {\rm Der} \nolimits}

\def\coker{\mathop {\rm coker} \nolimits}

\let\xrarr=\xrightarrow
\let\veps=\varepsilon
\let\rarr=\rightarrow

\let\mcal=\mathcal
\let\mfrak=\mathfrak
\let\eus=\EuScript

\newenvironment{enum}{\begin{list}{}{\topsep=2pt \itemsep=0pt
\parsep=0pt}}{\end{list}}

\title{Moduli spaces of flat Lie algebroid connections}

\author{Libor Křižka\thanks{The author of this article was supported by the grant LC 505 of the Ministry of Education, Youth and Sports (MŠMT).} \\
{\small Institute of Mathematics of the Academy of Sciences of the Czech Republic}  \\
{\small e-mail: {\tt krizka@karlin.mff.cuni.cz}}
}

\date{\today}

\begin{document}

\maketitle

\begin{abstract}
\noindent We shall prove that a moduli space of flat irreducible Lie algebroid connections over a compact manifold has locally a natural structure of a smooth differentiable space. This is a generalization of some well known results for the moduli space of holomorphic structures on a complex vector bundle over a compact complex manifold.


\end{abstract}

\thispagestyle{empty}

\section*{Introduction}

Moduli spaces arise naturally in classification problems in geometry. Typically, one has a set whose elements represent algebro-geometric objects of some fixed kind and an equivalence relation on this set saying when two such objects are the same in some sense, and the problem is to describe the set of equivalence classes. One would like to give the set of equivalence classes some structure of a geometric space (usually a smooth manifold, a scheme or an algebraic stack). If it can be done, then one can parametrize such objects by introducing coordinates on the resulting space.

The word \emph{moduli} is due to B.\,Riemann, who used it as a synonym for parameters when he showed that the space of equivalence classes of Riemann surfaces of a given genus $g$ (for $g > 1$) depends on $3g-3$ complex numbers. Moduli spaces were first understood as spaces of parameters rather than spaces of objects.

Moduli spaces have many important applications in mathematics and physics. In geometry, they allow us to construct invariants of manifolds, for example, Donaldson and Seiberg--Witten invariants of compact four-manifolds, Gromov--Witten invariants of symplectic manifolds and many others. Another application is the deformation theory describing small perturbations of a given object.

The motivation for the study of moduli spaces of flat Lie algebroid connections on vector bundles over compact manifolds is based on the two most important special cases. The first one is the moduli space of holomorphic structures on complex vector bundles over compact complex manifolds (more about this and the closely related Hitchin--Kobayashi correspondence can be found in \cite{Kobayashi1987}, \cite{Lubke_Okonek1987} and \cite{Lubke_Teleman1995}). The second one is the moduli space of Higgs bundles on compact Riemann surfaces, see \cite{Hitchin1987}, \cite{Simpson1988} and \cite{Simpson1992}.

The main purpose of this paper is to prove that certain moduli spaces of Lie algebroid connections on real (complex) vector bundles over compact manifolds carry a natural structure of real (complex) smooth differentiable spaces.

Let $L$ be a real (complex) Lie algebroid over a compact manifold $M$ such that $L$ satisfies the condition of ellipticity and let $E$ be a real (complex) vector bundle over $M$. In the first part of the paper, we develop the theory of Lie algebroid connections. In the second part, we give a differential geometric construction of the moduli space $\eus{M}^*\!(E,L)$ of isomorphism classes of irreducible flat Lie algebroid connections. This moduli space is a locally Hausdorff real (complex) smooth differentiable space of a finite dimension. We use the Kuranishi's method to describe local models for $\eus{M}^*\!(E,L)$. Locally around any point $[\nabla] \in \eus{M}^*\!(E,L)$ this space is a zero set of the smooth mapping
\begin{align*}
  \Phi \colon \mcal{O} \subset \mcal{H}^1(E,\nabla) \rarr \mcal{H}^2(E,\nabla)
\end{align*}
of finite dimensional cohomology spaces, and $\mcal{O}$ is an open subset in $\mcal{H}^1(E,\nabla)$.

\section{Lie algebroids -- definitions and examples}


Lie algebroids were first introduced and studied by J.\,Pradines \cite{Pradines1966a}, following the substantial work of C.\,Ehresmann and P.\,Libermann on \emph{differential groupoids} (later called \emph{Lie groupoids}), as infinitesimal objects for differential groupoids. Just as Lie algebras are the infinitesimal objects of Lie groups, Lie algebroids are the infinitesimal objects of Lie groupoids. They are generalization of both -- Lie algebras and tangent vector bundles.

\subsection{Lie algebroids}

\definition{A real (complex) Lie algebroid $(L \xrarr{\pi} M,[\cdot\,,\cdot],a)$ is a real (complex) vector bundle $\pi \colon L \rarr M$ together with a real (complex) Lie algebra bracket $[\cdot\,,\cdot]$ on the space of sections $\Gamma(M,L)$ and a homomorphism of vector bundles $a \colon L \rarr TM$ ($a \colon L \rarr TM_\mathbb{C}$), called the \emph{anchor map}, covering the identity on $M$, i.e., the following diagram
\begin{align*}
\bfig
\square[L`TM`M`M;a`\pi`\pi_M`{\rm id}_{\rm M}]
\efig
\quad \text{resp.}\quad
\bfig
\square[L`TM_\mathbb{C}`M`M;a`\pi`\pi_M`{\rm id}_M]
\efig
\end{align*}
commutes. Moreover, the anchor map fulfills
\begin{enum}
\item[i)] $a([\xi_1,\xi_2])=[a(\xi_1),a(\xi_2)]$ resp. $a([\xi_1,\xi_2])=[a(\xi_1),a(\xi_2)]_\mathbb{C}$,
\item[ii)] $[\xi_1,f\xi_2]=f[\xi_1,\xi_2]+(a(\xi_1)f)\xi_2$, (the Leibniz rule)
\end{enum}
for all $\xi_1,\xi_2 \in \Gamma(M,L)$ and $f \in C^\infty(M,\mathbb{R})$ resp. $f \in C^\infty(M,\mathbb{C})$.}

\definition{If $(L_1 \rarr M, [\cdot\,,\cdot]_{L_1},a_{L_1})$ and $(L_2 \rarr M, [\cdot\,,\cdot]_{L_2},a_{L_2})$ are Lie algebroids, then a vector bundle homomorphism $\varphi \colon L_1 \rarr L_2$ covering the identity on $M$ is a \emph{Lie algebroid morphism} if $a_{L_2} \circ \varphi = a_{L_1}$ and $\varphi$ induces a Lie algebra homomorphism form $\mfrak{X}_{L_1}(M)$ to $\mfrak{X}_{L_2}(M)$.}

\remark Let $A$ be a commutative $\mathbb{K}$-algebra\footnote{The letter $\mathbb{K}$ stands for the field $\mathbb{R}$ or $\mathbb{C}$.} with unit. We denote by $\Der_\mathbb{K}(A)$ the $A$-module of $\mathbb{K}$-linear derivations of $A$. Recall that $\Der_\mathbb{K}(A)$ is naturally a Lie algebra over $\mathbb{K}$ with respect to the usual commutator.

A Lie--Rinehart $A$-algebra is an $A$-module $L$ endowed with a structure of a Lie algebra over $\mathbb{K}$ and with a morphism $a \colon L \rarr \Der_\mathbb{K}(A)$ of $A$-modules, called the \emph{anchor map}, satisfying the following axioms:
\begin{enum}
  \item[i)] $a([x,y]_L)=[a(x),a(y)]$ for $x,y \in L$, i.e., $a$ is a morphism of Lie algebras over $\mathbb{K}$,
  \item[ii)] $[x,fy]_L=f[x,y]_L+(a(x)f)y$ for $x,y \in L$ and $f \in A$.
\end{enum}
\medskip

Consider the commutative $\mathbb{R}$-algebra $A=C^\infty(M,\mathbb{R})$, then $\Der_\mathbb{R}(A)$ is the Lie algebra of vector fields on $M$. Afterwards, the space of sections $\Gamma(M,L)$ of a real Lie algebroid $(L \rarr M,[\cdot\,,\cdot]_L,a)$ is a Lie--Rinehart $A$-algebra.

In fact, Lie--Rinehart algebras are the algebraic counterparts of Lie algebroids, just as modules over a ring are the algebraic counterpart of vector bundles.

\subsection{Examples of Lie algebroids}

Let us present now a few basic examples of Lie algebroids.

\subsubsection{Real Lie algebroids}

\example(tangent bundles) One of the trivial examples of a Lie algebroid over $M$ is the tangent bundle $L=TM$ of $M$, with the identity mapping as the anchor map and the Lie bracket of vector fields as the Lie bracket.
\medskip

\example(Lie algebras) Any real Lie algebra $\mfrak{g}$ is a real Lie algebroid over a one-point manifold, with zero anchor map.
\medskip

\example(foliations) Let $L \subset TM$ be an involutive regular distribution on a manifold $M$. Then $L$ has a Lie algebroid structure with the inclusion as the anchor map and the Lie bracket is the usual Lie bracket of vector fields. By the Frobenius theorem, the distribution $L$ gives a regular foliation on $M$. On the other hand, to any regular foliation on $M$ is associated an involutive regular distribution and therefore a Lie algebroid over $M$.
\medskip

\example(bundles of Lie algebras) A bundle of Lie algebras is a vector bundle $L \rarr M$ with a skew-symmetric $C^\infty(M,\mathbb{R})$-bilinear mapping $[\cdot\,,\cdot] \colon \Gamma(M,L) \times \Gamma(M,L) \rarr \Gamma(M,L)$, i.e., $[\cdot\,,\cdot] \in \Gamma(M,\Lambda^2L^* \otimes L)$, satisfying the Jacobi identity. If we define the anchor map by $a(\xi)=0$ for $\xi \in \Gamma(M,L)$, then $(L\rarr M, [\cdot\,,\cdot],a)$ is a Lie algebroid. On the other hand, any Lie algebroid with zero anchor map is a bundle of Lie algebras. Because $[\xi_1,f\xi_2]=f[\xi_1,\xi_2]+(a(\xi_1)f)\xi_2= f[\xi_1,\xi_2]$, we obtain $[\cdot\,,\cdot] \in \Gamma(M,\Lambda^2L^* \otimes L)$.

Note that the notion of a bundle of Lie algebras is weaker then of a Lie algebra bundle, when one requires that $L$ is locally trivial bundle of Lie algebras (in particular, all Lie algebras $L_x$ are isomorphic).
\medskip

\example(vector fields) Lie algebroid structures on the trivial real line bundle over $M$ are in a one-to-one correspondence with vector fields on $M$. Given a vector field $X \in \mfrak{X}(M)$, we denote by $L_X$ the induced Lie algebroid. As a vector bundle $L_X=M\times \mathbb{R}$. Because $\Gamma(M,L_X) \simeq C^\infty(M,\mathbb{R})$, the anchor map is given by the multiplication by $X$, i.e., $a(f)=fX$, and the Lie bracket of two sections $f,g \in \Gamma(M,L_X)$ is defined by
\begin{align}
  [f,g]=f\mcal{L}_X(g)-g \mcal{L}_X(f).
\end{align}

\example(action Lie algebroids) Consider an infinitesimal right action of a real Lie algebra $\mfrak{g}$ on a manifold $M$, i.e., a Lie algebra homomorphism $\zeta \colon \mfrak{g} \rarr \mfrak{X}(M)$. The usual situation is when we have a right action $r \colon M \times G \rarr M$ of a Lie group $G$ with the Lie algebra $\mfrak{g}$. Then
\begin{align}
  \zeta_X(x)=T_er_x.X= {{\rm d} \over {\rm d}t}_{|0}\, x.\exp(tX),
\end{align}
where $X\in \mfrak{g}$ and $x \in M$, gives an infinitesimal right action of $\mfrak{g}$ on $M$. We define a Lie algebroid $\mfrak{g} \ltimes M$, called the action Lie algebroid or the transformation Lie algebroid, by the following way. As a vector bundle $\mfrak{g} \ltimes M = M \times \mfrak{g}$, it is a trivial vector bundle over $M$. Seeing that $\Gamma(M,\mfrak{g} \ltimes M) \simeq C^\infty(M,\mfrak{g})$, the anchor map is given by
\begin{align}
  a(f)(x)=\zeta_{f(x)}(x),
\end{align}
while the Lie bracket on sections is defined by
\begin{align}
  [f,g](x)=[f(x),g(x)]_\mfrak{g} + (\zeta_{f(x)}g)(x)- (\zeta_{g(x)}f)(x).
\end{align}
The Lie bracket is uniquely determined by the Leibniz rule and the condition that
\begin{align}
  [c_X,c_Y]=c_{[X,Y]}
\end{align}
for all $X,Y \in \mfrak{g}$, where $c_X$ denotes the constant section of $M \times \mfrak{g}$.
\medskip

\example(two forms) For any closed 2-form $\omega \in \Omega^2(M,\mathbb{R})$, we define a Lie algebroid $L_\omega$ as follows. As a vector bundle $L_\omega = TM \oplus (M \times \mathbb{R})$, the anchor map is the projection on the first component, while the Lie bracket on sections $\Gamma(M,L_\omega)\simeq \mfrak{X}(M) \oplus C^\infty(M,\mathbb{R})$ is given by
\begin{align}
  [(X,f),(Y,g)]=([X,Y],\mcal{L}_X(g)- \mcal{L}_Y(f)+ \omega(X,Y)).
\end{align}

\example(Atiyah sequences) In 1957, Atiyah \cite{Atiyah1957} constructed in the setting of vector bundles for any principal fiber bundle $(P,p,M,G)$ the following key example of a transitive Lie algebroid $\mcal{A}(P)$, called the \emph{Atiyah algebroid}.

If $r \colon P \times G \rarr P$ is the principal right action, then $\hat{r} \colon TP \times G \rarr TP$ denotes the right action given by $\hat{r}^g=Tr^g$. The space of orbits $TP/G$ of the right action $\hat{r}$ has a canonical structure of a vector bundle over $M$. Moreover, we have an isomorphism $\Phi \colon \Gamma(M,TP/G) \xrarr{\sim} \mfrak{X}(P)^G$ of $C^\infty(M,\mathbb{R})$-modules, where $f\xi=(f \circ p)\xi$ for $f\in C^\infty(M,\mathbb{R})$ and $\xi \in \mfrak{X}(P)^G$. Furthermore, since $Tp$ is constant on orbits of the right action $\hat{r}$, the diagram
\begin{align*}
  \bfig
  \square[TP`TM`\mcal{A}(P)`TM;Tp`q`{\rm id}_{TM}`p_*]
  \efig
\end{align*}
commutes for a uniquely determined smooth mapping $p_* \colon \mcal{A}(P) \rarr TM$, where $q\colon TP \rarr \mcal{A}(P)$ is the quotient map.

As a vector bundle $\mcal{A}(P)=TP/G$, the anchor map is $p_*$, while the Lie bracket on sections $\Gamma(M,\mcal{A}(P))$ is given by
\begin{align}
  [\xi_1,\xi_2]=\Phi^{-1}([\Phi(\xi_1),\Phi(\xi_2)])
\end{align}
for $\xi_1,\xi_2 \in \Gamma(M,\mcal{A}(P))$.

Because $\mcal{A}(P)$ is a transitive Lie algebroid, i.e., the anchor map $p_*$ is surjective, we obtain the following short exact sequence
\begin{align}
  0 \xrarr{\phantom{p_*}} \ad(P) \xrarr{\phantom{p_*}} \mcal{A}(P) \xrarr{p_*} TM \xrarr{\phantom{p_*}} 0
\end{align}
of Lie algebroids over $M$ known as the Atiyah sequence associated to the principal $G$-bundle $P$, where the Lie bracket on $\Gamma(M,\ad(P))$ is induced from the given one on $\Gamma(M,\mcal{A}(P))$.

If $L$ is a transitive Lie algebroid over $M$, then the associated short exact sequence
\begin{align}
  0 \xrarr{\phantom{a}} \ker a \xrarr{i}  L \xrarr{a} TM \xrarr{\phantom{a}} 0
  \label{abstract atiyah sequence}
\end{align}
of Lie algebroids is called the \emph{abstract Atiyah sequence}. Note that not all abstract Atiyah sequences come from sequences associated to principal fiber bundles.
\medskip

\example(Poisson manifolds) Any Poisson structure on a manifold $M$ induces, in a natural way, a Lie algebroid structure on the cotangent bundle $T^*\!M$ of $M$. Let $\pi \in \Gamma(M,\Lambda^2TM)$ be a Poisson bivector on $M$, which is related to the Poisson bracket by $\{f,g\}=\pi(df, dg)$.
If we use the notation
\begin{align}
  \pi^\sharp \colon T^*\!M \rarr TM
\end{align}
for the mapping defined by $\beta(\pi^\sharp(\alpha))=\pi(\alpha,\beta)$ for $\alpha,\beta \in \Omega^1(M,\mathbb{R})$, then the Hamiltonian vector field $X_f$ associated to a smooth function $f$ on $M$ is defined by $X_f=\pi^\sharp(df)$. The anchor map is $\pi^\sharp$ and the Lie bracket is given by
\begin{align}
  [\alpha,\beta]=\mcal{L}_{\pi^\sharp(\alpha)}(\beta) - \mcal{L}_{\pi^\sharp(\beta)} (\alpha)- d(\pi(\alpha,\beta)).
\end{align}
This Lie algebroid structure on $T^*\!M$ is the unique one with the property that $a(df)=X_f$ and $[df,dg]=d\{f,g\}$ for all $f,g \in C^\infty(M,\mathbb{R})$.
When $\pi$ is nondegenerate, $M$ is a symplectic manifold and this Lie algebra structure of $\Gamma(M,T^*\!M)$ is isomorphic to that of $\Gamma(M,TM)$.
\medskip

\example(Nijenhuis manifolds) Let $M$ be a manifold with a Nijenhuis structure, i.e., a vector valued 1-form $\mcal{N} \in \Omega^1(M,TM)$ with the vanishing Nijenhuis torsion. Recall that the Nijenhuis torsion $T_\mcal{N}\in \Omega^2(M,TM)$ is defined by
\begin{align}
  T_\mcal{N}(X,Y)=[\mcal{N}X,\mcal{N}Y]-\mcal{N}[\mcal{N}X,Y]- \mcal{N}[X,\mcal{N}Y] + \mcal{N}^2[X,Y]
\end{align}
for $X,Y \in \mfrak{X}(M)$, note that $T_\mcal{N}={1\over 2} [\mcal{N},\mcal{N}]$ for the Frölicher--Nijenhuis bracket. A vector valued 1-form $\mcal{N}$ is called a \emph{Nijenhuis tensor} if its Nijenhuis torsion vanishes. To any Nijenhuis structure $\mcal{N}$, there is associated a new Lie algebroid structure on $TM$. The anchor map is given by $a(X)=\mcal{N}(X)$, while the Lie bracket is defined by
\begin{align}
  [X,Y]_\mcal{N}=[\mcal{N}X,Y]+[X,\mcal{N}Y]-\mcal{N}[X,Y].
\end{align}
It is well known that powers of Nijenhuis tensors, considered as endomorphisms of the tangent bundle, are Nijenhuis tensors. Also any complex structure $\mcal{J}$ on $M$ is a Nijenhuis tensor.
\medskip

\example(generalized Nijenhuis manifolds) Let $(L \rarr M,[\cdot\,,\cdot],a)$ be a Lie algebroid and let $\mcal{N} \colon L \rarr L$ be a homomorphism of vector bundles covering the identity on $M$, such that its Nijenhuis torsion vanishes, i.e.,
\begin{align}
  [\mcal{N}X,\mcal{N}Y]-\mcal{N}[\mcal{N}X,Y]- \mcal{N}[X,\mcal{N}Y] + \mcal{N}^2[X,Y]=0
\end{align}
for all $X,Y \in \Gamma(M,L)$. If we define the anchor map by $a_\mcal{N}(X)=(a \circ \mcal{N})(X)$ and the Lie bracket by
\begin{align}
  [X,Y]_\mcal{N}=[\mcal{N}X,Y]+[X,\mcal{N}Y]-\mcal{N}[X,Y],
\end{align}
then this gives a new Lie algebroid structure on $L$.
\medskip

\example(trivial Lie algebroids) For any real Lie algebra $\mfrak{g}$, we define a Lie algebroid $L_\mfrak{g}$ over a manifold $M$ by the following way. As a vector bundle $L_\mfrak{g}=TM \oplus (M \times \mfrak{g})$, the anchor map is the projection on the first component and the Lie bracket on sections $\Gamma(M,L_\mfrak{g}) \simeq \mfrak{X}(M) \oplus C^\infty(M,\mfrak{g})$ is defined by
\begin{align}
  [(X,f),(Y,g)]=([X,Y],[f,g]),
\end{align}
where the bracket on sections $\Gamma(M,M\times \mfrak{g}) \simeq C^\infty(M,\mfrak{g})$ is given by
\begin{align}
  [f,g](x)=[f(x),g(x)]_\mfrak{g}.
\end{align}


\example(jet prolongation of Lie algebroids) Let $(L \xrarr{\pi} M,[\cdot\,,\cdot],a)$ be a Lie algebroid, then the $r$-th jet prolongation $J^r\!L$ of $L$ for $r\in \mathbb{N}_0$ has a unique Lie algebroid structure. The anchor map is given by $a_{J^r\!L}= \pi^r_0 \circ a$, where $\pi^r_0 \colon J^r\!L \rarr L$ is the canonical projection, while the Lie bracket is uniquely determined by requiring that the $r$-th jet prolongation
\begin{align}
  j^r \colon \Gamma(M,L) \rarr \Gamma(M,J^r\!L)
\end{align}
is a homomorphism of Lie algebroids. More about the relation of jet prolongation of Lie algebroids to the Cartan's method of equivalence one can find in \cite{Blaom2005}.

\subsubsection{Complex Lie algebroids}

\example(Lie algebras) Any complex Lie algebra $\mfrak{g}$ is a complex Lie algebroid over a one-point manifold, with zero anchor map.
\medskip

\example(complexification of real Lie algebroids) Let $(L \rarr M,[\cdot\,,\cdot],a)$ be a real Lie algebroid, then $L_\mathbb{C}$ becomes a complex Lie algebroid when the Lie bracket is extended complex bilinearly and the anchor map is extended complex linearly.
\medskip

\example(complex manifolds) Let $M$ be a manifold with a complex structure, i.e., a vector valued $1$-form $\mcal{J} \in \Omega^1(M,TM)$ with the vanishing Nijenhuis torsion and satisfying $\mcal{J}^2=-{\rm id}_{TM}$. As a vector bundle $L=TM^{0,1}$, the anchor map $a \colon TM^{0,1} \rarr TM_\mathbb{C}= TM^{1,0} \oplus TM^{0,1}$ is inclusion only and the Lie bracket is the Lie bracket of complexified vector fields.
\medskip

\example(CR manifolds) An abstract CR structure on a manifold $M$ is a complex subbundle $L$ of the complexified tangent bundle $TM_\mathbb{C}$ satisfying  $L \cap \bar{L}= 0$ and $[\xi_1,\xi_2]_\mathbb{C} \in \Gamma(M,L)$ for any $\xi_1,\xi_2 \in \Gamma(M,L)$. The bundle $L$ is called an abstract CR structure on the manifold $M$. The CR codimension of the CR structure is $k=\dim M - 2\rk L$. In the case $k=1$, the CR structure is said to be hypersurface type. The complex vector bundle $L$ together with the Lie bracket of vector fields and the inclusion of $L$ into $TM_\mathbb{C}$ as the anchor map is a complex Lie algebroid.
\medskip

\example(involutive structures) The previous two examples are special cases of more general involutive structures. An involutive structure on a manifold $M$ is a complex subbundle $L$ of $TM_\mathbb{C}$ such that $[\xi_1,\xi_2]_\mathbb{C} \in \Gamma(M,L)$ for any $\xi_1,\xi_2 \in \Gamma(M,L)$. The complex vector bundle $L$ together with the Lie bracket of vector fields and the inclusion of $L$ into $TM_\mathbb{C}$ as the anchor map is a complex Lie algebroid.
\medskip

\section{Differential geometry of Lie algebroids}

Because we can think of Lie algebroids as generalized tangent bundles, we may use similar constructions for them.

Given a real (complex) Lie algebroid $(L \xrarr{\pi} M,[\cdot\,,\cdot],a)$. A section of the vector bundle $\Lambda^k L^*$ for $k\in \mathbb{N}_0$ is called a $k$-form of $L$ and the space of all $k$-forms will be denoted by $\Omega^k_L(M)$. Similarly, a section of the vector bundle $\Lambda^k L$ for $k\in \mathbb{N}_0$ is called a $k$-vector field of $L$ and the space of all $k$-vector fields will be denoted by $\mfrak{X}^k_L(M)$. Let $\Omega^k_L(M)=\{0\}$ and $\mfrak{X}^k_L(M)=\{0\}$ for $k < 0$, then we denote by
\begin{align}
  \Omega^\bullet_L(M)=\bigoplus_{k\in \mathbb{Z}}\, \Omega^k_L(M) \quad \text{and} \quad \mfrak{X}^\bullet_L(M)= \bigoplus_{k\in \mathbb{Z}}\, \mfrak{X}^k_L(M)
\end{align}
the graded vector space of forms of $L$ and the graded vector space of multivector fields of $L$, respectively. For a real (complex) vector bundle $E \rarr M$ a section of the vector bundle $\Lambda^kL^* \otimes E$ is called $E$-valued $k$-form of $L$. The space of sections will be denoted by $\Omega^k_L(M,E)$.

The graded vector space $\Omega^\bullet_L(M)$ has a natural structure of a graded commutative algebra via the wedge product
\begin{align}
  (\omega \wedge \tau)(\xi_1,\dots,\xi_{p+q})={1 \over p!\,q!} \sum_\sigma \sign(\sigma)\cdot\omega(\xi_{\sigma(1)},\dots,\xi_{\sigma(p)})\, \tau(\xi_{\sigma(p+1)},\dots,\xi_{\sigma(p+q)}),
\end{align}
where $\omega \in \Omega^p_L(M)$, $\tau \in \Omega^q_L(M)$ and $\xi_1,\dots,\xi_{p+q} \in \mfrak{X}_L(M)$.

Further, there is a differential operator $d_L \colon \Omega^\bullet_L(M) \rarr \Omega^{\bullet+1}_L(M)$ on the graded commutative algebra $\Omega^\bullet_L(M)$ defined by
\begin{multline}
  (d_L \omega)(\xi_0,\dots,\xi_k)=\sum_{i=0}^k (-1)^i a(\xi_i)\omega(\xi_0,\dots, \widehat{\xi_i}, \dots,\xi_k) \\
   + \sum_{0\leq i<j \leq k} (-1)^{i+j}\omega([\xi_i,\xi_j],\xi_0,\dots, \widehat{\xi_i}, \dots, \widehat{\xi_j},\dots,\xi_k)
\end{multline}
for $\omega \in \Omega^k_L(M)$ and $\xi_0,\dots,\xi_k \in \mfrak{X}_L(M)$.
The linear differential operator $d_L$ is called the \emph{Lie algebroid differential} of $L$ or simply the \emph{de Rham differential} of $L$.
Besides for any $\xi \in \mfrak{X}_L(M)$ we define the \emph{insertion operator} $i^L_\xi  \colon \Omega^\bullet_L(M) \rarr \Omega^{\bullet-1}_L(M)$ by
\begin{align}
(i^L_\xi \omega)(\xi_1,\dots,\xi_k)=\omega(\xi,\xi_1,\dots,\xi_k)
\end{align}
and the \emph{Lie derivative} $\mcal{L}^L_\xi \colon \Omega^\bullet_L(M) \rarr \Omega^\bullet_L(M)$ through
\begin{align}
  (\mcal{L}^L_\xi \omega)(\xi_1,\dots,\xi_k)= a(\xi) \omega(\xi_1,\dots,\xi_k) - \sum_{i=1}^k \omega(\xi_1,\dots, [\xi,\xi_i], \dots, \xi_k)
\end{align}
for $\omega \in \Omega^k_L(M)$ and $\xi,\xi_1,\dots,\xi_k \in \mfrak{X}_L(M)$.
\medskip

\remark As $\Omega^\bullet_L(M)$ is a graded commutative algebra, the graded vector space of graded derivations
\begin{align}
  \Der\Omega^\bullet_L(M) = \bigoplus_{k \in \mathbb{Z}} \Der_k\!\Omega^\bullet_L(M),
\end{align}
where $\Der_k \!\Omega^\bullet_L(M)$ is the vector space of graded derivations of degree $k$, has a structure of a graded Lie algebra with the Lie bracket  defined by
\begin{align}
  [D_1,D_2]= D_1 \circ D_2 - (-1)^{k_1 k_2} D_2 \circ D_1
\end{align}
for $D_1 \in \Der_{k_1}\! \Omega^\bullet_L(M)$ and $D_2 \in \Der_{k_2}\! \Omega^\bullet_L(M)$.
\medskip

\lemma{The insertion operator $i^L_\xi \colon \Omega^\bullet_L(M) \rarr \Omega^{\bullet-1}_L(M)$ and the Lie derivative $\mcal{L}^L_\xi \colon \Omega^\bullet_L(M) \rarr \Omega^\bullet_L(M)$ have the following properties:
\begin{enum}
  \item[i)] $i^L_\xi(\omega \wedge \tau)= i^L_\xi\omega \wedge \tau + (-1)^{\deg(\omega)} \omega \wedge i^L_\xi\tau$, i.e., $i^L_\xi$ is a graded derivation od degree $-1$,
  \item[ii)] $\mcal{L}^L_\xi(\omega \wedge \tau)= \mcal{L}^L_\xi\omega \wedge \tau +  \omega \wedge \mcal{L}^L_\xi\tau$, i.e., $\mcal{L}^L_\xi$ is a graded derivation od degree $0$,
  \item[iii)] $[\mcal{L}^L_\xi, i^L_\eta]=i^L_{[\xi,\eta]}$,
  \item[iv)] $[\mcal{L}^L_\xi,\mcal{L}^L_\eta]= \mcal{L}^L_{[\xi,\eta]}$,
  \item[v)]  $[i^L_\xi, i^L_\eta]=0$.
\end{enum}}

\proof{The proof goes along the same lines as in the special case $L=TM$, see \cite{Kolar_Michor_Slovak1993}.}

\lemma{The Lie algebroid differential $d_L \colon \Omega^\bullet_L(M) \rarr \Omega^{\bullet+1}_L(M)$ has the following properties:
\begin{enum}
  \item[i)] $d_L(\omega \wedge \tau)= d_L\omega \wedge \tau + (-1)^{\deg(\omega)} \omega \wedge d_L\tau$, i.e., $d_L$ is a graded derivation od degree $1$,
  \item[ii)] $d_L \circ d_L ={1 \over 2}[d_L,d_L]=0$, i.e., $d_L$ is a differential,
  \item[iii)] $[\mcal{L}^L_\xi,d]=0$,
  \item[iv)] $[i^L_\xi, d]=\mcal{L}^L_\xi$ (Cartan's formula).
\end{enum}}

\proof{The proof goes along the same lines as in the special case $L=TM$, see \cite{Kolar_Michor_Slovak1993}.}

Because $d_L$ is a graded derivation of degree $1$ and a differential, i.e., $d^2_L=0$, the graded commutative algebra $\Omega^\bullet_L(M)$ is a differential graded commutative algebra. The cohomology of the complex
\begin{align}
  0 \xrarr{\phantom{d_L}} \Omega^0_L(M) \xrarr{d_L} \Omega^1_L(M) \xrarr{d_L} \dots \xrarr{d_L} \Omega^r_L(M) \xrarr{\phantom{d_L}} 0, \label{de Rham complex}
\end{align}
where $r = \rk L$, called the \emph{Lie algebroid cohomology} of $L$, we will  denote it by $H^\bullet_L(M)$. It unifies de Rham and Chevalley--Eilenberg cohomology. When $L=TM$, we obtain $H^\bullet_{TM}(M)= H^\bullet_{{\rm dR}}(M)$, on the other hand, when $L=\mfrak{g}$, i.e., $L$ is a Lie algebroid over a one-point manifold, we receive $H^\bullet_\mfrak{g}(M)=H^\bullet(\mfrak{g},\mfrak{g})$.
Moreover, because $d_L$ is a graded derivation of degree $1$, the Lie algebroid cohomology $H^\bullet_L(M)$ of $L$ is a graded commutative $\mathbb{R}$-algebra ($\mathbb{C}$-algebra).

Furthermore, we can ask when is this complex an elliptic complex? For any $f \in C^\infty(M,\mathbb{R})$ and $\omega \in \Omega^k_L(M)$ we have
\begin{align*}
(\ad(f)d_L)\, \omega = d_L(f\omega)-fd_L \omega = d_Lf \wedge \omega + fd_L \omega -fd_L \omega =a^*\!(df) \wedge \omega.
\end{align*}
Hence for the principal symbol $\sigma_1(d_L)$ we get
\begin{align*}
\sigma_1(d_L)(\xi_x)=a^*\!(\xi_x)\,\wedge  \colon (\Lambda^kL^*)_x \rarr (\Lambda^{k+1}L^*)_x
\end{align*}
for every $x\in M$ and $\xi_x \in T^*_x\!M$, i.e., the symbol is the exterior multiplication by $a^*\!(\xi_x)$. Therefore, we obtain the \emph{Koszul complex}
\begin{align*}
  0 \xrightarrow{\phantom{a^*\!(\xi_x)\wedge}} (\Lambda^0L^*)_x \xrightarrow{a^*\!(\xi_x)\wedge} (\Lambda^1L^*)_x \xrightarrow{a^*\!(\xi_x)\wedge} \hdots \xrightarrow{a^*\!(\xi_x)\wedge} (\Lambda^rL^*)_x \xrightarrow{\phantom{a^*\!(\xi_x)\wedge}} 0,
\end{align*}
where $r=\rk L$, which is an exact sequence if and only if $a^*\!(\xi_x) \neq 0$.
Thus, the differential complex \eqref{de Rham complex} is elliptic if and only if the corresponding Koszul complex is an exact sequence for any $x \in M$ and $0 \neq \xi_x \in T^*_x\!M$, in other words if and only if $a^*\!(\xi_x) \neq 0$ for any $x \in M$ and $0 \neq \xi_x \in T^*_x\!M$.

If $L \xrarr{a} TM$ is a real Lie algebroid, then the elipticity is equivalent to the requirement that $a^* \colon T^*\!M \rarr L^*$ is injective or that $a \colon L \rarr TM$ is surjective.
For a complex Lie algebroid $L \xrarr{a} TM_\mathbb{C}$ it corresponds to the requirement that ${a^*}_{|T^*\!M} \colon T^*\!M \hookrightarrow (TM_\mathbb{C})^* \rarr L^*$ is injective.
\medskip

\lemma{\label{lemma-kernel}Let $M$ be a compact manifold and let $L \rarr M$ be a real (complex) Lie algebroid satisfying the ellipticity condition. Then $\ker d_L$ consists of locally constant real (complex) functions.}

\proof{Because the Lie algebroid $L$ satisfies the ellipticity condition, the complex
\begin{align*}
  0 \xrarr{\phantom{d_L}} \Omega^0_L(M) \xrarr{d_L} \Omega^1_L(M) \xrarr{d_L} \dots \xrarr{d_L} \Omega^r_L(M) \xrarr{\phantom{d_L}} 0,
\end{align*}
where $r = \rk L$, is an elliptic complex. Further, since $M$ is a compact manifold, $\dim \ker d_L < \infty$. As we can write $d_L(fg)=g(d_Lf) + f(d_Lg)$, we get that $\ker d_L \subset C^\infty(M,\mathbb{K})$ is a finite dimensional $\mathbb{K}$-algebra.

Suppose that there exists $f \in C^\infty(M,\mathbb{K})$ such that $d_Lf=0$ and $f$ is not locally constant function on $M$. But the subset $\{1,f,f^2,\dots,f^n\}$ of $\ker d_L$ is linearly independent for any $n \in \mathbb{N}$ which is a contradiction with the fact that $\ker d_L$ is a finite dimensional vector space.}

\section{Linear Lie algebroid connections}

In this section we introduce the notion of linear Lie algebroid connections, i.e., Lie algebroid connections on vector bundles. It is a natural generalization of linear connections on vector bundles. Since Lie algebroids can be understood as generalized tangent bundles, it is possible to use analogous constructions for linear Lie algebroid connections as for linear connections.
\medskip

\remark The letter $\mathbb{K}$ stands for the field $\mathbb{R}$ of real numbers or the field $\mathbb{C}$ of complex numbers.
\medskip

\definition{Let $(L \rarr M,[\cdot\,,\cdot],a)$ be a real (complex) Lie algebroid and let $E \rarr M$ be a real (complex) vector bundle. We will denote the space of sections of the vector bundle $\Lambda^kL^* \otimes E$ for $k \in \mathbb{N}_0$ by $\Omega_L^k(M,E)$ and  sections will be called $E$-valued $k$-forms of $L$ or $k$-forms of $L$ with values in $E$. A \emph{linear Lie algebroid connection for $L$} or an \emph{$L$-connection} on the vector bundle $E$ is a $\mathbb{K}$-linear mapping
\begin{align}
  \nabla \colon \Omega_L^0(M,E) \rarr \Omega_L^1(M,E)
\end{align}
satisfying Leibniz rule $\nabla(fs)=d_Lf \otimes s + f\nabla s$ for any $f \in C^\infty(M,\mathbb{K})$ and $s \in \Omega_L^0(M,E)$.}

\remark For any $\xi \in \mfrak{X}_L(M)$ we have a $\mathbb{K}$-linear mapping $\nabla_{\!\xi} \colon \Omega^0_L(M,E) \rarr \Omega^0_L(M,E)$ given by
\begin{align}
  \nabla_{\!\xi}s= i^L_\xi(\nabla s)
\end{align}
for $s \in \Omega^0_L(M,E)$, called the \emph{covariant derivative} along $\xi$. Moreover, it satisfies
\begin{gather}
  \nabla_{\!\xi}(fs)=(\mcal{L}^L_\xi \!f)s + f \nabla_{\!\xi}s, \qquad  \nabla_{\!f\xi}s=f \nabla_{\!\xi}s \label{covariant derivative - Leibniz} \\
  \intertext{and}
  \nabla_{\!\xi_1+\xi_2}s=\nabla_{\!\xi_1}s + \nabla_{\!\xi_2}s  \label{covariant derivative - linearity}
\end{gather}
for all $f \in C^\infty(M,\mathbb{K})$, $\xi,\xi_1,\xi_2 \in \mfrak{X}_L(M)$ and $s \in \Omega^0_L(M,E)$.
Therefore, a linear Lie algebroid connection for $L$ on the vector bundle $E$ can be equivalently defined as a $\mathbb{K}$-bilinear mapping
\begin{gather}
  \nabla \colon \mfrak{X}_L(M) \times \Omega^0_L(M,E) \rarr \Omega^0_L(M,E), \nonumber \\
  (\xi,s) \mapsto \nabla_{\!\xi}s
\end{gather}
satisfying \eqref{covariant derivative - Leibniz} for all $\xi \in \mfrak{X}_L(M)$, $f \in C^\infty(M,\mathbb{K})$ and $s \in \Omega^0_L(M,E)$.
\medskip

Tensorial operations on vector bundles may be extended naturally to vector bundles with $L$-connections. More precisely, if $E_1$ and $E_2$ are two vector bundles with $L$-connections $\nabla^{E_1}$ and $\nabla^{E_2}$, then $E_1 \otimes E_2$ has a naturally induced $L$-connection $\nabla^{E_1 \otimes E_2}$ uniquely determined by the formula
\begin{equation}
  \nabla^{E_1 \otimes E_2}_{\!\xi}(s_1 \otimes s_2)=\nabla^{E_1}_{\!\xi}s_1 \otimes s_2 + s_1 \otimes \nabla^{E_2}_{\!\xi} s_2
\end{equation}
for all $\xi \in \mfrak{X}_L(M)$, $s_1 \in \Omega^0_L(M,E_1)$ and $s_2 \in \Omega^0_L(M,E_2)$. If we are given a vector bundle $E$ with an $L$-connection $\nabla^E$, then the dual vector bundle $E^*$ has a natural $L$-connection $\nabla^{E^*}$ defined by the identity
\begin{equation}
  \mcal{L}_\xi^L \!\langle t,s \rangle = \langle \nabla_{\!\xi}^{E^*}\! t,s \rangle + \langle t, \nabla^E_{\!\xi} s \rangle
\end{equation}
for all $\xi \in \mfrak{X}_L(M)$, $s \in \Omega^0_L(M,E)$ and $t \in \Omega^0_L(M,E^*)$, where $\langle \cdot\,,\cdot \rangle \colon \Omega^0_L(M,E^*) \times \Omega^0_L(M,E) \rarr C^\infty(M,\mathbb{K})$ is the natural pairing. In particular, any $L$-connection $\nabla^E$ on a vector bundle $E$ induces an $L$-connection $\nabla^{\End(E)}$ on $\End(E) \simeq E^* \otimes E$ by the rule
\begin{equation}
  (\nabla^{\End(E)}_{\!\xi}T)s=\nabla^E_{\!\xi}(Ts)-T(\nabla^E_{\!\xi} s)=[\nabla^E_{\!\xi},T]s
\end{equation}
for all $\xi \in \mfrak{X}_L(M)$, $T \in \Omega^0_L(M,\End(E))$ and $s \in \Omega^0_L(M,E)$.
\medskip

For any vector bundle $E$ the graded vector space $\Omega^\bullet_L(M,E)$ is a graded $\Omega^\bullet_L(M)$-module through
\begin{align}
  (\alpha \wedge \omega)(\xi_1,\dots,\xi_{p+q})={1 \over p!\,q!} \sum_{\sigma} \sign(\sigma)\cdot \alpha(\xi_{\sigma(1)},\dots,\xi_{\sigma(p)})\, \omega(\xi_{\sigma(p+1)},\dots,\xi_{\sigma(p+q)}), \label{multiplication-module}
\end{align}
where $\alpha \in \Omega^p_L(M)$, $\omega \in \Omega^q_L(M,E)$ and $\xi_1,\dots,\xi_{p+q} \in \mfrak{X}_L(M)$. The graded module homomorphisms $\Phi \colon \Omega^\bullet_L(M,E) \rarr \Omega^\bullet_L(M,E)$ (so that $\Phi(\alpha \wedge \omega)= (-1)^{\deg(\Phi) \deg(\alpha)}\, \alpha \wedge \Phi(\omega)$) coincide with the mappings $\mu(A)$ for $A \in \Omega^p_L(M,\End(E))$, which are given by
\begin{align}
 \big(\mu(A)\,\omega\big)(\xi_1,\dots,\xi_{p+q})={1 \over p!\,q!} \sum_{\sigma} \sign (\sigma) \cdot A(\xi_{\sigma(1)},\dots,\xi_{\sigma(p)})\,\omega(\xi_{\sigma(p+1)}, \dots, \xi_{\sigma(p+q)}),
\end{align}
where $\omega \in \Omega^q_L(M,E)$ and $\xi_1,\dots,\xi_{p+q} \in \mfrak{X}_L(M)$.
Further, the graded vector space $\Omega^\bullet_L(M,\End(E))$ has a natural structure of a graded associative algebra via
\begin{align}
  (\omega \wedge \tau)(\xi_1,\dots,\xi_{p+q})={1 \over p!\,q!} \sum_{\sigma} \sign(\sigma)\cdot (\omega(\xi_{\sigma(1)},\dots,\xi_{\sigma(p)}) \circ \tau(\xi_{\sigma(p+1)},\dots,\xi_{\sigma(p+q)})) \label{multiplication}
\end{align}
and a natural structure of a graded Lie algebra through
\begin{align}
  [\omega, \tau](\xi_1,\dots,\xi_{p+q})={1 \over p!\,q!} \sum_{\sigma} \sign(\sigma)\cdot [\omega(\xi_{\sigma(1)},\dots,\xi_{\sigma(p)}), \tau(\xi_{\sigma(p+1)},\dots,\xi_{\sigma(p+q)})], \label{lie-bracket}
\end{align}
where $\omega \in \Omega^p_L(M,\End(E))$, $\tau \in \Omega^q_L(M,\End(E))$ and $\xi_1,\dots,\xi_{p+q} \in \mfrak{X}_L(M)$. Comparing these two definitions we may write
\begin{align}
  [\omega,\tau]= \omega \wedge \tau - (-1)^{\deg(\omega) \deg(\tau)} \,\tau \wedge \omega.
\end{align}
for $\omega, \tau \in \Omega^\bullet_L(M,\End(E))$.
\medskip


Let $\nabla$ be an $L$-connection on a vector bundle $E$, then the \emph{covariant exterior derivative}
\begin{align}
d^\nabla \colon \Omega_L^\bullet(M,E) \rarr \Omega_L^{\bullet+1}(M,E)
\end{align}
is defined by
\begin{multline}
 (d^\nabla\! \omega)(\xi_0,\xi_1,\dots,\xi_k)= \sum_{i=0}^k (-1)^i \nabla_{\!\xi_i} \omega(\xi_0,\dots,\hat{\xi_i},\dots,\xi_k) \\
  + \sum_{0\leq i<j \leq k} (-1)^{i+j} \omega([\xi_i,\xi_j],\xi_0,\dots,\hat{\xi_i},\dots,\hat{\xi_j},\dots,\xi_k),
\end{multline}
where $\omega \in \Omega^k_L(M,E)$ and $\xi_0,\dots,\xi_k \in \mfrak{X}_L(M)$.
\medskip

\lemma{The covariant exterior derivative $d^\nabla \colon \Omega_L^\bullet(M,E) \rarr \Omega_L^{\bullet+1}(M,E)$ has the following properties:
\begin{enum}
  \item[i)] $d^\nabla\! (\Omega^k_L(M,E)) \subset \Omega^{k+1}_L(M,E)$,
  \item[ii)] ${d^\nabla}_{|\Omega^0_L(M,E)}=\nabla$,
  \item[iii)] $d^\nabla\!(\alpha \wedge \omega)= d_L\alpha \wedge \omega + (-1)^{\deg(\alpha)} \alpha \wedge d^\nabla\! \omega$ for $\alpha \in \Omega^\bullet_L(M)$
      and $\omega \in \Omega^\bullet_L(M,E)$ (the graded Leibniz rule),
  \item[iv)] $d^{\nabla^{\End(E)}}[\omega,\tau]= [d^{\nabla^{\End(E)}}\omega,\tau] + (-1)^{\deg(\omega)}[\omega, d^{\nabla^{\End(E)}}\tau]$ for $\omega, \tau \in \Omega^\bullet_L(M,\End(E))$.
\end{enum}}

\proof{Properties i) and ii) follow immediately from the definition.

\noindent
iii) It suffices to investigate decomposable forms $\omega= \beta \otimes s$ for $\beta \in \Omega^q_L(M)$ and $s \in \Omega^0_L(M,E)$. From the definition we obtain  $d^\nabla\!(\beta \otimes s)=d_L\beta \otimes s + (-1)^q \beta \wedge d^\nabla\!s$. Afterwards for $\alpha \in \Omega^p_L(M)$ we have
\begin{align*}
  d^\nabla\!(\alpha \wedge (\beta \otimes s))&= d^\nabla\!((\alpha \wedge \beta) \otimes  s)= d_L(\alpha \wedge \beta) \otimes s + (-1)^{p+q}(\alpha \wedge \beta) \wedge d^\nabla\!s \\
  &= (d_L\alpha \wedge \beta) \otimes s + (-1)^p(\alpha \wedge d_L\beta) \otimes s + (-1)^{p+q}(\alpha \wedge \beta) \wedge d^\nabla\!s \\
  &=d_L\alpha \wedge (\beta \otimes s) + (-1)^p\alpha \wedge d^\nabla\!(\beta \otimes s).
\end{align*}

\noindent
iv) For decomposable forms $\omega=\alpha \otimes s$, $\tau= \beta \otimes t$, where $s,t \in \Omega^0_L(M,\End(E))$, $\alpha \in \Omega^p_L(M)$ and $\beta \in \Omega^q_L(M)$, we have
$[\alpha \otimes s, \beta \otimes t]=(\alpha \wedge \beta) \otimes [s,t]$. Hence we can write
\begin{align*}
  d^{\nabla^{\End(E)}}[\alpha \otimes s, \beta \otimes t]&= d^{\nabla^{\End(E)}} ((\alpha \wedge \beta) \otimes [s,t]) \\
  &=d_L(\alpha \wedge \beta) \otimes [s,t] + (-1)^{p+q} (\alpha \wedge \beta ) \wedge d^{\nabla^{\End(E)}}[s,t] \\
  &=(d_L\alpha \wedge \beta) \otimes [s,t] + (-1)^p (\alpha \wedge d_L\beta) \otimes [s,t] \\
  &\quad + (-1)^{p+q} (\alpha \wedge \beta ) \wedge [d^{\nabla^{\End(E)}}s,t] + (-1)^{p+q} (\alpha \wedge \beta ) \wedge [s,d^{\nabla^{\End(E)}}t] \\
  &=[d_L\alpha \otimes s, \beta \otimes t] + (-1)^p[\alpha \otimes s, d_L\beta \otimes t] + (-1)^p [\alpha \wedge d^{\nabla^{\End(E)}}s, \beta \otimes t] \\
  &\quad + (-1)^{p+q}[\alpha \otimes s,\beta \wedge d^{\nabla^{\End(E)}} t] \\
  &=[d^{\nabla^{\End(E)}}(\alpha \otimes s), \beta \otimes t]+ (-1)^p [\alpha \otimes s, d^{\nabla^{\End(E)}}(\beta \otimes t)],
\end{align*}
where we used that $d^{\nabla^{\End(E)}}[s,t]=[d^{\nabla^{\End(E)}}s,t] + [s,d^{\nabla^{\End(E)}}t]$ which follows from the classical Jacobi identity for $\mathbb{K}$-linear mappings on $\Omega^0_L(M,E)$, thus we are done.}

\lemma{Denote by $\eus{A}(E,L)$ the set of all $L$-connections on a vector bundle $E$. Then $\eus{A}(E,L)$ is an affine space modeled on the vector space $\Omega^1_L(M,\End(E))$.}

\proof{We first prove that $\eus{A}(E,L)$ is non-empty. Because on the vector bundle $E$ there exists a connection $\tilde{\nabla} \colon \Omega^0(M,E) \rarr \Omega^1(M,E)$, we may define an $L$-connection $\nabla \colon \Omega^0_L(M,E) \rarr \Omega^1_L(M,E)$ by the formula
\begin{align*}
  \nabla_{\!\xi}s=\tilde{\nabla}_{\!a(\xi)}s
\end{align*}
for $\xi \in \mfrak{X}_L(M)$ and $s \in \Omega^0_L(M,E)$. The rest of the proof is very simple. We need to verify that, if $\nabla$ and $\nabla'$ are two $L$-connections, then $(\nabla'-\nabla) \colon \Omega^0_L(M,E) \rarr \Omega^1_L(M,E)$ is a $C^\infty(M,\mathbb{K})$-linear mapping.
As $(\nabla' -\nabla)(fs)=d_Lf \otimes s + f \nabla'\! s - d_Lf \otimes s - f \nabla s = f(\nabla'-\nabla)s$, there exists a uniquely determined $\alpha \in \Omega^1_L(M,\End(E))$ such that $\nabla' -\nabla = \mu(\alpha)$.}

\remark Thus, if we fix some $\nabla_{\!0}$ in $\eus{A}(E,L)$, we may write
\begin{align}
  \eus{A}(E,L)=\{\nabla_{\!0} + \mu(\alpha);\, \alpha \in \Omega^1_L(M,\End(E))\}. \label{konexe}
\end{align}
This description will permit us to define Sobolev completions of $\eus{A}(E,L)$.
\medskip

\definition{If we are given an $L$-connection $\nabla$ on a vector bundle $E$, then the \emph{curvature} $R^\nabla \in \Omega^2_L(M,\End(E))$ of the $L$-connection $\nabla$ is defined by the formula
\begin{equation}
R^\nabla\!(\xi,\eta)s=\nabla_{\!\xi}\nabla_{\!\eta}s -\nabla_{\!\eta}\nabla_{\!\xi}s - \nabla_{\![\xi,\eta]}s = [\nabla_{\!\xi},\nabla_{\!\eta}]s - \nabla_{\![\xi,\eta]}s,
\end{equation}
where $\xi,\eta \in \mfrak{X}_L(M)$ and $s \in \Omega^0_L(M,E)$.

An $L$-connection with zero curvature is called a \emph{flat $L$-connection}. We denote the set of all flat $L$-connections on a vector bundle $E$ by $\eus{H}(E,L)$.}

\lemma{\label{lemma-curvature}Let $\nabla$ be an $L$-connection on a vector bundle $E$, then
\begin{align}
  (d^\nabla \!\circ d^\nabla) \omega =\mu(R^\nabla)\omega
\end{align}
for all $\omega \in \Omega^\bullet_L(M,E)$.}

\proof{First we verify that $R^\nabla\!(\xi,\eta)s=(d^\nabla\!(d^\nabla \!s))(\xi,\eta)$. This is a consequence upon the following computation
\begin{align*}
\begin{split}
(d^\nabla\!(d^\nabla \!s))(\xi,\eta)&=\nabla_{\!\xi}((d^\nabla \!s)(\eta))- \nabla_{\!\eta}((d^\nabla \!s)(\xi)) - (d^\nabla \!s)([\xi,\eta]) \\
  &=\nabla_{\!\xi} \nabla_{\!\eta} s - \nabla_{\!\eta} \nabla_{\!\xi} s -\nabla_{\![\xi,\eta]}s \\
  &= R^\nabla\!(\xi,\eta)s
\end{split}
\end{align*}
for all $\xi,\eta \in \mfrak{X}_L(M)$ and $s \in \Omega^0_L(M,E)$.
Further, it suffices to investigate only decomposable forms $\omega =\alpha \otimes s$ for $\alpha \in \Omega^k_L(M)$ and $s \in \Omega^0_L(M,E)$. Afterwards, we can write
\begin{equation*}
\begin{split}
  (d^\nabla\! \circ d^\nabla) (\alpha \otimes s)&=d^\nabla\! (d_L \alpha \otimes s + (-1)^k \alpha \wedge d^\nabla\! s) \\
  &= 0 + (-1)^{k+1} d_L\alpha \wedge d^\nabla\!s + (-1)^k d_L \alpha \wedge d^\nabla \!s+ (-1)^{2k} \alpha \wedge (d^\nabla\! \circ d^\nabla)s \\
  &= \alpha \wedge \mu(R^\nabla)\,s \\
  &= \mu(R^\nabla)\,(\alpha \otimes s)
\end{split}
\end{equation*}
hence we have got  $d^\nabla\! \circ d^\nabla = \mu(R^\nabla)$ and this finishes the proof.}

Given an $L$-connection on a vector bundle $E$, the mapping $\nabla \colon \Omega^0_L(M,E) \rarr \Omega^1_L(M,E)$ can be extended to the following sequence of first order differential operators
\begin{align}
0 \xrightarrow{\phantom{d^\nabla}} \Omega^0_L(M,E) \xrightarrow{d^\nabla}
\Omega^1_L(M,E) \xrightarrow{d^\nabla} \hdots \xrightarrow{d^\nabla} \Omega^r_L(M,E) \xrightarrow{\phantom{d^\nabla}} 0, \label{complex2}
\end{align}
where $r=\rk L$. It is a differential complex if and only if the curvature $R^\nabla$ of the $L$-connection $\nabla$ is zero, i.e., $\nabla$ is a flat $L$-connection.

A natural question is when is this differential complex an elliptic complex? Let $f \in C^\infty(M,\mathbb{R})$, then we may write
\begin{align*}
(\ad(f)d^\nabla)\omega=d^\nabla\!(f\omega)- f d^\nabla\! \omega = d_Lf \wedge \omega + f d^\nabla\! \omega - f d^\nabla\! \omega= a^*\!(df) \wedge \omega
\end{align*}
for any $\omega \in \Omega^k_L(M,E)$. Therefore for the principal symbol $\sigma_1\!(d^\nabla)$ we obtain
\begin{align*}
\sigma_1(d^\nabla)(\xi_x)=a^*\!(\xi_x)\, \wedge \colon (\Lambda^kL^* \otimes E)_x \rarr (\Lambda^{k+1}L^* \otimes E)_x
\end{align*}
for every $x \in M$ and $\xi_x \in T^*_x\!M$, i.e., the symbol is the exterior multiplication by $a^*\!(\xi_x)$. Hence, we have the \emph{twisted Koszul complex}
\begin{align*}
  0 \xrightarrow{\phantom{a^*\!(\xi_x)\wedge}} (\Lambda^0L^*\otimes E)_x \xrightarrow{a^*\!(\xi_x)\wedge} \hdots \xrightarrow{a^*\!(\xi_x)\wedge} (\Lambda^rL^*\otimes E)_x \xrightarrow{\phantom{a^*\!(\xi_x)\wedge}} 0,
\end{align*}
where $r=\rk L$, which is an exact sequence if and only if $a^*\!(\xi_x) \neq 0$.
Thus, the differential complex \eqref{complex2} is elliptic if and only if the corresponding twisted Koszul complex is an exact sequence for any $x \in M$ and $0 \neq \xi_x \in T^*_x\!M$, in other words if and only if $a^*\!(\xi_x) \neq 0$ for any $x \in M$ and $0 \neq \xi_x \in T^*_x\!M$.

If $L \xrarr{a} TM$ is a real Lie algebroid, then the elipticity is equivalent to the requirement that $a^* \colon T^*\!M \rarr L^*$ is injective or that $a \colon L \rarr TM$ is surjective.
For a complex Lie algebroid $L \xrarr{a} TM_\mathbb{C}$ it corresponds to the requirement that ${a^*}_{|T^*\!M} \colon T^*\!M \hookrightarrow (TM_\mathbb{C})^* \rarr L^*$ is injective.
These are the same conditions as for the ellipticity of the complex \eqref{de Rham complex}. We will call this condition the \emph{ellipticity condition} for the Lie algebroid $L$.
\medskip

\lemma{If $\nabla$ is an $L$-connection on a vector bundle $E$, then we have
\begin{align}
  d^{\nabla^{\End(E)}} R^\nabla =0.
\end{align}
This is called the \emph{Bianchi identity} for $R^\nabla$.}

\proof{For any $\xi_1,\xi_2,\xi_3 \in \mfrak{X}_L(M)$ we may write
\begin{align*}
  (d^{\nabla^{\End(E)}} R^\nabla)(\xi_1,\xi_2,\xi_3) &=[\nabla_{\!\xi_1}, R^\nabla\!(\xi_2,\xi_3)] - [\nabla_{\!\xi_2}, R^\nabla\!(\xi_1,\xi_3)] + [\nabla_{\!\xi_3}, R^\nabla\!(\xi_1,\xi_2)] \\
  & \quad - R^\nabla\!([\xi_1,\xi_2],\xi_3) + R^\nabla\!([\xi_1,\xi_3],\xi_2) - R^\nabla\!([\xi_2,\xi_3],\xi_1) \\
  & =\sum_{cykl}\,([\nabla_{\!\xi_1},[\nabla_{\!\xi_2}, \nabla_{\!\xi_3}]] - [\nabla_{\!\xi_1}, \nabla_{\![\xi_2,\xi_3]}])
    - \sum_{cykl}\,([\nabla_{\![\xi_1,\xi_2]},\nabla_{\!\xi_3}] - \nabla_{\![[\xi_1,\xi_2],\xi_3]}) \\
  &=  -\sum_{cykl}\, [\nabla_{\!\xi_1}, \nabla_{\![\xi_2,\xi_3]}]
    - \sum_{cykl}\,[\nabla_{\![\xi_1,\xi_2]},\nabla_{\!\xi_3}]
  =0,
\end{align*}
where we used the classical Jacobi identity for commutators of $\mathbb{K}$-linear mappings.}

\lemma{\label{connection-curvature}Consider two $L$-connections $\nabla, \nabla'$ on a vector bundle $E$. Then
\begin{align}
  R^{\nabla'}&=R^\nabla+ d^{\nabla^{\End(E)}}\alpha + \alpha \wedge \alpha \label{connection-curvature1}\\
  &=R^\nabla+ d^{\nabla^{\End(E)}}\alpha + {1 \over 2}\, [\alpha, \alpha], \label{connection-curvature2}
\end{align}
where $\alpha \in \Omega^1_L(M,\End(E))$ such that  $\nabla'-\nabla=\mu(\alpha)$.}

\proof{The proof is a straightforward computation only. We have
\begin{equation*}
  \begin{split}
    R^{\nabla'}\!\!(\xi,\eta)&=[\nabla'_{\!\xi}, \nabla'_{\!\eta}] - \nabla'_{\![\xi,\eta]} \\
    &=[\nabla_{\!\xi}+\alpha(\xi), \nabla_{\!\eta}+ \alpha(\eta)] - (\nabla_{\![\xi,\eta]}+ \alpha([\xi,\eta])) \\
    &= [\nabla_{\!\xi}, \nabla_{\!\eta}] - \nabla_{\![\xi,\eta]}
    + [\nabla_{\!\xi},\alpha(\eta)] - [\nabla_{\!\eta}, \alpha(\xi)] - \alpha([\xi,\eta]) + [\alpha(\xi),\alpha(\eta)] \\
    &=R^\nabla\!(\xi,\eta) + \nabla^{\End(E)}_{\!\xi}\alpha(\eta)- \nabla^{\End(E)}_{\!\eta}\alpha(\xi) - \alpha([\xi,\eta]) + [\alpha(\xi),\alpha(\eta)] \\
    &=R^\nabla\!(\xi,\eta) + (d^{\nabla^{\End(E)}}\alpha)(\xi,\eta)  + (\alpha \wedge \alpha)(\xi,\eta) \\
    &=R^\nabla\!(\xi,\eta) + (d^{\nabla^{\End(E)}}\alpha)(\xi,\eta)  + \textstyle{{1 \over 2}}[\alpha,\alpha](\xi,\eta)
    \end{split}\end{equation*}
for all $\xi, \eta \in \mfrak{X}_L(M)$, so we are done.}
\medskip

Therefore, if we fix some flat $L$-connection $\nabla_{\!0} \in \eus{H}(E,L)$, then, using the result of Lemma \ref{connection-curvature}, we may write
\begin{align}
  \eus{H}(E,L)=\{\nabla_{\!0} + \mu(\alpha);\, \alpha \in \Omega^1_L(M,\End(E)),\,  d^{\nabla_{\!0}^{\End(E)}} \alpha + \alpha \wedge \alpha =0\}.
\end{align}
This description, similarly like in the case of $\eus{A}(E,L)$, will allow us to define Sobolev completions of $\eus{H}(E,L)$.

\section{Group of gauge transformations}

Let $E \xrarr{p} M$ be a real (complex) vector bundle, then a vector bundle homomorphism is a smooth mapping $\varphi \colon E \rarr E$ such that there exists mapping $\underline{\varphi} \colon M \rarr M$, the diagram
\begin{align*}
 \bfig
 \square[E`E`M`M;\varphi`p`p`\underline{\varphi}]
 \efig
\end{align*}
commutes and for each $x \in M$ the mapping $\varphi_x = \varphi_{|E_x} \colon E_x \rarr E_{\underline{\varphi}(x)}$ is $\mathbb{K}$-linear. As $p \colon E \rarr M$ is a fibered manifold and $\underline{\varphi} \circ p$ is smooth, we get that $\underline{\varphi}$ is smooth. Moreover, if $\varphi \colon E \rarr E$ is a vector bundle automorphism, then the previous diagram commutes for a uniquely determined diffeomorphism $\underline{\varphi} \colon M \rarr M$.

In case we denote by $\Aut(E)$ the group of vector bundle automorphisms and by $\Diff(M)$ the group of diffeomorphisms of $M$, then  we get the group homomorphism from $\Aut(E)$ into $\Diff(M)$ defined by $\varphi \mapsto \underline{\varphi}$. The kernel $\Gau(E)$ of this homomorphism is called the \emph{group of gauge transformations} and its elements are called \emph{gauge transformations}. Thus $\Gau(E)$ is the group of vector bundle automorphisms $\varphi \colon E \rarr E$ satisfying $p \circ \varphi = p$. Hence, we have the following exact sequence
\begin{align}
  \{e\} \rarr \Gau(E) \rarr \Aut(E) \rarr \Diff(M)
\end{align}
of groups.

Furthermore, we define the \emph{Lie algebra of gauge transformations} $\gau(E)$. As a vector space it is $\Omega^0_L(M,\End(E))$, while the Lie bracket is given by
\begin{align}
  [\gamma_1,\gamma_2]=\gamma_1 \circ \gamma_2 - \gamma_2 \circ \gamma_1
\end{align}
for $\gamma_1, \gamma_2 \in \Omega^0_L(M,\End(E))$.
\medskip

The group of gauge transformations $\Gau(E)$ has a natural left action on the vector space $\Omega^k_L(M,\End(E))$ defined through
\begin{align}
  (\Ad_\varphi(\omega))(\xi_1,\dots,\xi_k)=\varphi \circ \omega(\xi_1,\dots, \xi_k) \circ \varphi^{-1}, \label{Ad action}
\end{align}
where $\varphi \in \Gau(E)$, $\omega \in \Omega^k_L(M,\End(E))$ and $\xi_1,\dots,\xi_k \in \mfrak{X}_L(M)$. Further, this gives a left action of the Lie algebra of gauge transformations $\gau(E)$ on $\Omega^k_L(M,\End(E))$ via
\begin{align}
  \ad_\gamma(\omega)=[\gamma,\omega]
\end{align}
for $\gamma \in \mfrak{gau}(E)$ and $\omega \in \Omega^k_L(M,\End(E))$. So we have got a representation of $\Gau(E)$ and $\mfrak{gau}(E)$ on the graded vector space $\Omega^\bullet_L(M,\End(E))$.
\medskip

\remark Furthermore, there is a left action of the group $\Aut(E)$ on the space of sections $\Omega^0_L(M,E)$ defined by
\begin{align}
  \varphi \cdot s = \varphi \circ s \circ \underline{\varphi}^{-1},
\end{align}
where $\varphi \in \Aut(E)$ and $s \in \Omega^0_L(M,E)$.

\section{Change of connections}

Let $(L \rarr M,[\cdot\,,\cdot],a)$ be a real (complex) Lie algebroid and let $E \rarr M$ be a real (complex) vector bundle. Further, consider a gauge transformation $\varphi$ and an $L$-connection $\nabla$ on $E$. We can define a $\mathbb{K}$-bilinear mapping $\nabla^\varphi \colon \mfrak{X}_L(M) \times \Omega^0_L(M,E) \rarr \Omega^0_L(M,E)$ by
\begin{align}
  \nabla^\varphi_{\!\xi}\!s= \varphi^{-1} \circ \nabla_{\!\xi}(\varphi \circ s)
\end{align}
for any $\xi \in \mfrak{X}_L(M)$ and $s \in \Omega^0_L(M,E)$. Since we may write
\begin{align*}
   \nabla^\varphi_{\!\xi}\!(fs)&=\varphi^{-1}\circ \nabla_{\!\xi}(\varphi \circ (fs))= \varphi^{-1}\circ \nabla_{\!\xi}(f(\varphi\circ s)) \\
  &= \varphi^{-1}\circ ((\mcal{L}^L_\xi \!f)(\varphi\circ s) +f\nabla_{\!\xi}(\varphi\circ s)) \\
  &= (\mcal{L}^L_\xi \!f)s  +f(\varphi^{-1} \circ \nabla_{\!\xi}(\varphi\circ s)) \\
  &= (\mcal{L}^L_\xi \!f)s  +f\nabla^\varphi_{\!\xi}\!s
\end{align*}
and moreover we have
\begin{align*}
\nabla^\varphi_{\!f\xi}s= \varphi^{-1}\circ \nabla_{\!f\xi}(\varphi\circ s) = \varphi^{-1}\circ (f\nabla_{\!\xi}(\varphi\circ s)) = f(\varphi^{-1}\circ \nabla_{\!\xi}(\varphi\circ s)) = f \nabla^\varphi_{\!\xi}\!s
\end{align*}
for all $\xi \in \mfrak{X}_L(M)$, $f \in C^\infty(M,\mathbb{K})$ and $s \in \Omega^0_L(M,E)$, therefore $\nabla^\varphi$ is an $L$-connection on $E$.

As $\nabla^\varphi$ is an $L$-connection, we can define a natural right action of $\Gau(E)$ on the space $\eus{A}(E,L)$ of $L$-connections by
\begin{align}
  (\nabla,\varphi) \mapsto  \nabla \cdot \varphi =\nabla^\varphi. \label{gauge group right action}
\end{align}
It is easy to see that this really defines a right action.
\medskip

\lemma{\label{lemma curvature-transformation}Let $\nabla$ be an $L$-connection on $E$, then we have
\begin{align}
  R^{\nabla^\varphi}=\Ad_{\varphi^{-1}}(R^\nabla) \label{curvature-transformation}
\end{align}
for any gauge transformation $\varphi \in \Gau(E)$.}

\proof{It follows immediately that
\begin{align*}
  \begin{split}
    R^{\nabla^\varphi}\!\!(\xi,\eta)&=[\nabla^\varphi_{\!\xi}, \nabla^\varphi_{\!\eta}]- \nabla^\varphi_{\![\xi,\eta]} \\
    & = \varphi^{-1} \circ [\nabla_{\!\xi},\nabla_{\!\eta}] \circ \varphi - \varphi^{-1} \circ \nabla_{\![\xi,\eta]} \circ \varphi \\
    & = \varphi^{-1} \circ R^\nabla\!(\xi,\eta) \circ \varphi
  \end{split}
\end{align*}
for all $\xi,\eta \in \mfrak{X}_L(M)$.}

Because $\eus{H}(E,L)$ is invariant under the action of $\Gau(E)$, as it follows from Lemma~\ref{lemma curvature-transformation}, we have a right action of $\Gau(E)$ on the space of flat $L$-connections $\eus{H}(E,L)$. Therefore, we define the \emph{moduli space}
\begin{align}
  \eus{B}(E,L) =  \lower-0.5ex\hbox{$\eus{A}(E,L)$}\mkern-2mu \big/ \mkern-2mu \lower0.5ex\hbox{$\Gau(E)$} \label{moduli-konexe}
\end{align}
of gauge equivalence classes of $L$-connections and the \emph{moduli space}
\begin{align}
  \eus{M}(E,L) =  \lower-0.5ex\hbox{$\eus{H}(E,L)$}\mkern-2mu \big/ \mkern-2mu \lower0.5ex\hbox{$\Gau(E)$} \label{moduli-flat-konexe}
\end{align}
of gauge equivalence classes of flat $L$-connections.
\medskip

Now we take up the question of reducible $L$-connections. Given an $L$-connection $\nabla$ on a vector bundle $E$, then the \emph{isotropy subgroup} or the \emph{stabilizer subgroup} of $\nabla$ is the subgroup $\Gau(E)_\nabla$ of $\Gau(E)$ that leaves $\nabla$ fixed, i.e.,
\begin{align}
\Gau(E)_\nabla  =\{\varphi \in \Gau(E);\, \nabla \cdot \varphi = \nabla\}.
\end{align}
Every such group contains the subgroup $\Gau(E)^0$ defined by
\begin{align*}
\Gau(E)^0=\{f\cdot {\rm id}_E;\, f\in C^\infty(M,\mathbb{K}),\, d_Lf=0,\, f(x) \neq 0\ \forall\, x \in M\}.
\end{align*}
Moreover, it is a normal subgroup of the group of gauge transformations $\Gau(E)$. Provided that $M$ is a connected manifold and $L$ is a transitive Lie algebroid, then $\Gau(E)^0= \mathbb{K}^*\! \cdot {\rm id}_E$. Furthermore, we define the Lie algebra $\gau(E)^0$. As a vector space it is
\begin{align*}
  \gau(E)^0=\{f \cdot {\rm id}_E;\, f\in C^\infty(M,\mathbb{K}),\, d_Lf=0\},
\end{align*}
while the Lie bracket is trivial, i.e., $\gau(E)^0$ is an abelian Lie algebra.
\medskip

\definition{An $L$-connection $\nabla$ on a vector bundle $E$ is called \emph{irreducible} or \emph{simple}, in case $\Gau(E)_\nabla =\Gau(E)^0$, otherwise $\nabla$ is called \emph{reducible}. We denote the set of irreducible (simple) $L$-connections by $ \eus{A}^*\!(E,L) $ and the set of irreducible flat $L$-connections by $\eus{H}^*\!(E,L)$.}

\lemma{\label{lemma-irred}Let $\nabla$ be an $L$-connection on a vector bundle $E$ over a compact manifold $M$. Then the following conditions are equivalent:
\begin{enum}
  \item[i)] $\Gau(E)_\nabla= \Gau(E)^0$,
  \item[ii)] $\ker \nabla^{\End(E)}= \gau(E)^0$,
\end{enum}
where $\nabla^{\End(E)}$ is the induced $L$-connection on $\End(E)$.}

\proof{Consider a gauge transformation $\varphi \in \Gau(E)$. Then the requirement $\nabla \cdot \varphi= \nabla$ means that for any $\xi \in \mfrak{X}_L(M)$ we have $\varphi^{-1} \circ \nabla_{\!\xi} \circ \varphi= \nabla_{\!\xi}$ and this is equivalent to $[\nabla_{\!\xi},\varphi]=0$. Therefore, we have got that $\varphi \in \Gau(E)_\nabla$ if and only if $\nabla^{\End(E)} \varphi =0$ and $\varphi \in \Gau(E)$.

Suppose that $\varphi \in \Gau(E)_\nabla$, then $\nabla^{\End(E)} \varphi =0$ and, provided that the condition ii) holds, we obtain $\varphi = f\cdot{\rm id}_E$ for some $f \in C^\infty(M,\mathbb{K})$ satisfying $d_Lf=0$. Hence we get $\Gau(E)_\nabla \subset \Gau(E)^0$ and because the converse inclusion is trivial, we have proved ii) $\Rightarrow$ i).

To prove the opposite implication, we use the compactness of the manifold $M$. Assume that $\varphi \in \ker \nabla^{\End(E)}$. Because $M$ is compact, there exists $c \in \mathbb{K}$ (with $|c|$ sufficiently large) so that $c \cdot {\rm id}_E + \varphi \in \Gau(E)$. Moreover, we have $\nabla^{\End(E)}(c \cdot {\rm id}_E + \varphi)=0$ and from the previous consideration, it follows $c \cdot {\rm id}_E + \varphi \in \Gau(E)_\nabla$. Besides, if we suppose $\Gau(E)_\nabla = \Gau(E)^0$, we obtain $\ker \nabla^{\End(E)} \subset \gau(E)^0$. The converse inclusion is trivial.}

From the fact that $\Gau(E)_{{\nabla^\varphi}}= \varphi^{-1} \cdot \Gau(E)_\nabla \cdot \varphi$ for any $\varphi \in \Gau(E)$ and $\nabla \in \eus{A}(E,L)$, we obtain that $\eus{A}^*\!(E,L)$ is invariant under the action of $\Gau(E)$ and the same for $\eus{H}^*\!(E,L)$. Thus, we can define, similarly as in \eqref{moduli-konexe} and \eqref{moduli-flat-konexe}, the moduli space
\begin{align}
  \eus{B}^*\!(E,L) =  \lower-0.5ex\hbox{$\eus{A}^*\!(E,L)$}\mkern-2mu \big/ \mkern-2mu \lower0.5ex\hbox{$\Gau(E)$} \label{moduli-ir-konexe}
\end{align}
of gauge equivalence classes of irreducible $L$-connections and the moduli space
\begin{align}
  \eus{M}^*\!(E,L) =  \lower-0.5ex\hbox{$\eus{H}^*\!(E,L)$}\mkern-2mu \big/ \mkern-2mu \lower0.5ex\hbox{$\Gau(E)$} \label{moduli-ir-flat-konexe}
\end{align}
of gauge equivalence classes of irreducible flat $L$-connections.

Because $\Gau(E)^0$ is a normal subgroup of $\Gau(E)$, we define the \emph{reduced group of gauge transformations} $\Gau(E)^{\rm r}$ by
\begin{align}
  \Gau(E)^{\rm r}= \lower-0.5ex\hbox{$\Gau(E)$}\mkern-2mu \big/ \mkern-2mu \lower0.5ex\hbox{$\Gau(E)^0$}. \label{reduced-gauge-group}
\end{align}
Then the right action of $\Gau(E)$ on $\eus{A}(E,L)$ factors trough the action of the reduced group of gauge transformations $\Gau(E)^{\rm r}$ since the group $\Gau(E)^0$ acts trivially on $\eus{A}(E,L)$, similarly for $\eus{H}(E,L)$. Therefore for the moduli spaces \eqref{moduli-konexe}, \eqref{moduli-flat-konexe} of $L$-connections we may write
\begin{align}
    \eus{B}(E,L) =  \lower-0.5ex\hbox{$\eus{A}(E,L)$}\mkern-2mu \big/ \mkern-2mu \lower0.5ex\hbox{$\Gau(E)^{\rm r}$} \qquad \text{and} \qquad
    \eus{M}(E,L) =  \lower-0.5ex\hbox{$\eus{H}(E,L)$}\mkern-2mu \big/ \mkern-2mu \lower0.5ex\hbox{$\Gau(E)^{\rm r}$}
\end{align}
and similarly  for the moduli spaces \eqref{moduli-ir-konexe}, \eqref{moduli-ir-flat-konexe} of irreducible $L$-connections we have
\begin{align}
    \eus{B}^*\!(E,L) =  \lower-0.5ex\hbox{$\eus{A}^*\!(E,L)$}\mkern-2mu \big/ \mkern-2mu \lower0.5ex\hbox{$\Gau(E)^{\rm r}$}
    \qquad \text{and} \qquad
    \eus{M}^*\!(E,L) =  \lower-0.5ex\hbox{$\eus{H}^*\!(E,L)$}\mkern-2mu \big/ \mkern-2mu \lower0.5ex\hbox{$\Gau(E)^{\rm r}$}.
\end{align}
The set $\eus{A}^*\!(E,L)$ of irreducible $L$-connections is the maximal subset of $\eus{A}(E,L)$ on which the reduced group of gauge transformations $\Gau(E)^{\rm r}$ acts freely, likewise for $\eus{H}^*\!(E,L)$.
\medskip

If we are given a gauge transformation $\varphi \in \Gau(E)$ and an $L$-connection $\nabla$ on a vector bundle $E$, then for the changed $L$-connection $\nabla^\varphi$ we have
\begin{align}
  \nabla^\varphi_{\!\xi}&= \nabla_{\!\xi} + \varphi^{-1} \circ \nabla^{\End(E)}_{\!\xi} \varphi
 =\nabla_{\!\xi}- \nabla^{\End(E)}_{\!\xi} \varphi^{-1} \circ \varphi,
\end{align}
where $\xi \in \mfrak{X}_L(M)$. The last equality follows by differentiating the identity $\varphi^{-1} \circ \varphi = {\rm id}_E$. More generally, if we fix some $L$-connection $\nabla$ and express another $L$-connection $\nabla'$ as $\nabla' =\nabla + \mu(\alpha)$, then
\begin{align}
\nabla'^\varphi_{\!\xi} = \nabla_{\!\xi} + \varphi^{-1} \circ \nabla^{\End(E)}_{\!\xi} \varphi + \varphi^{-1} \circ \alpha(\xi) \circ \varphi, \label{action}
\end{align}
hence, writing $\nabla'^\varphi = \nabla + \mu(\alpha^\varphi)$, we obtain
\begin{align}
  \alpha^\varphi\!(\xi)=\varphi^{-1} \circ \nabla^{\End(E)}_{\!\xi} \varphi + \varphi^{-1} \circ \alpha(\xi) \circ \varphi
\end{align}
for $\xi \in \mfrak{X}_L(M)$. This can be rewritten as
\begin{align}
  \alpha^\varphi &= \varphi^{-1} \wedge \nabla^{\End(E)} \varphi + \Ad_{\varphi^{-1}}(\alpha)
  \\ &= -\nabla^{\End(E)} \varphi^{-1} \wedge \varphi + \Ad_{\varphi^{-1}}(\alpha)
\end{align}
for $\varphi \in \Gau(E)$.

\section{Moduli spaces}

The moduli spaces $\eus{B}(E,L)$, $\eus{M}(E,L)$, $\eus{B}^*\!(E,L)$ and $\eus{M}^*\!(E,L)$ introduced in the previous section were only sets of gauge equivalence classes of $L$-connections. The main goal of this paper is to specify a geometric structure on these sets.

\medskip
From now on, we will assume that $M$ is a connected compact manifold. To endow the sets of gauge equivalence classes of $L$-connections with some geometric structure it is most convenient, and standard practise, to work within the framework of Sobolev spaces.

Let $(L \rarr M, [\cdot\,,\cdot],a)$ be a real (complex) Lie algebroid satisfying the ellipticity condition and let $E \rarr M$ be a real (complex) vector bundle. Further, consider a Riemannian metric $g$ on $M$ and denote by $h_E$ and $h_L$ an Euclidean (Hermitian) metric on $E$ and $L$, respectively. These metrics induce natural metrics on $E^*$, $\End(E) \simeq E^* \otimes E$, $\Lambda^k L^* \otimes \End(E)$ and others. The metric $g$ on $M$ defines the density $\vol(g)$ of the Riemannian metric and even induces a (regular) Borel measure $\mu_g$ on $M$. Therefore, we can construct appropriate Sobolev completions $L^{\ell,p}(M,\Lambda^k L^* \otimes \End(E))$ of $\Gamma(M,\Lambda^k L^* \otimes\End(E))$ defined for $1 \leq p < + \infty$ and $\ell \in \mathbb{N}_0$.
The corresponding Hilbert spaces for $p=2$ will be denoted by $\Omega^k_L(M,\End(E))_\ell$.
\medskip

For a fixed $L$-connection $\nabla_{\!0}$ on a vector bundle $E$, we define Sobolev completions $\eus{A}(E,L)_\ell$ of the space of $L$-connections $\eus{A}(E,L)$ for $\ell \in \mathbb{N}_0$, using \eqref{konexe}, as
\begin{align}
  \eus{A}(E,L)_\ell=\{\nabla_{\!0}+ \alpha;\, \alpha \in \Omega^1_L(M,\End(E))_\ell\}.
\end{align}
As $\eus{A}(E,L)_\ell$ is an affine Hilbert space, it is in fact a Hilbert manifold modeled on $\Omega^1_L(M,\End(E))_\ell$ whose tangent space at a point $\nabla$ is
\begin{align}
  T_\nabla \eus{A}(E,L)_\ell = \Omega^1_L(M,\End(E))_\ell.
\end{align}
Sobolev completions of the group of gauge transformations $\Gau(E)$ take a bit more work since it can not be identified with the space of sections of any vector bundle, nevertheless we have $\Gau(E) \subset \Omega^0_L(M,\End(E))$. If $\ell +1 > {1 \over 2} \dim M$, the Sobolev space $\Omega^0_L(M,\End(E))_{\ell+1}$ consists of continuous sections and, using the Sobolev multiplication theorem, we obtain that the product $\varphi \cdot \psi=\varphi \circ \psi$ in $\Omega^0_L(M,\End(E))$ can be extended to a continuous bilinear mapping
\begin{align}
  \Omega^0_L(M,\End(E))_{\ell+1} \times \Omega^0_L(M,\End(E))_{\ell+1} \rarr \Omega^0_L(M,\End(E))_{\ell+1}.
\end{align}
Because the set of all invertible elements is an open subset of $\Omega^0_L(M,\End(E))_{\ell+1}$ and  forms a topological group under multiplication, we define $\Gau(E)_{\ell+1}$ by
\begin{align}
\Gau(E)_{\ell +1}=\{\varphi \in \Omega^0_L(M,\End(E))_{\ell+1};\, \exists\, \psi \in \Omega^0_L(M,\End(E))_{\ell+1},\, \varphi \cdot \psi= \psi \cdot \varphi = {\rm id}_E \}.
\end{align}
Since $\Gau(E)_{\ell+1}$ is an open subset of the Hilbert space $\Omega^0_L(M,\End(E))_{\ell+1}$, the group $\Gau(E)_{\ell+1}$ is a Hilbert manifold modeled on $\Omega^0_L(M,\End(E))_{\ell+1}$. In fact, it is easy to see that $\Gau(E)_{\ell+1}$ is a Hilbert--Lie group with the Lie algebra
\begin{align}
  \gau(E)_{\ell+1}=\Omega^0_L(M,\End(E))_{\ell+1},
\end{align}
where the Lie bracket is given by
\begin{align}
  [\gamma_1,\gamma_2]=\gamma_1 \cdot \gamma_2 - \gamma_2 \cdot \gamma_1
\end{align}
for all $\gamma_1, \gamma_2  \in \Omega^0_L(M,\End(E))_{\ell+1}$.


The multiplication on the graded vector space $\Omega^\bullet_L(M,\End(E))$ defined by \eqref{multiplication} is extended, using the Sobolev multiplication theorem, to a continuous bilinear mapping on the graded Hilbert space $\Omega^\bullet_L(M,\End(E))_k$ in the range $k > {1 \over 2} \dim M$. With this bilinear mapping
\begin{gather}
  \Omega^p_L(M,\End(E))_k \times \Omega^q_L(M,\End(E))_k \rarr \Omega^{p+q}_L(M,\End(E))_k, \nonumber \\
  (\varphi, \psi ) \mapsto \varphi \wedge \psi,
\end{gather}
the graded vector space $\Omega^\bullet_L(M,\End(E))_k$ is a graded associative algebra.

Using the formula \eqref{action}, we extend the action of $\Gau(E)$ on $\eus{A}(E,L)$ to the action of $\Gau(E)_{\ell+1}$ on $\eus{A}(E,L)_\ell$ defined via
\begin{align}
 \nabla \cdot \varphi  = (\nabla_{\!0} + \alpha) \cdot \varphi &=\nabla_{\!0} + \varphi^{-1} \wedge (d^{\nabla_{\!0}} \varphi) + \Ad_{\varphi^{-1}}(\alpha),
\end{align}
where $\alpha \in \Omega^1_L(M,\End(E))_\ell$, $d^{\nabla_{\!0}} \colon \Omega^0_L(M,\End(E))_{\ell+1} \rarr \Omega^1_L(M,\End(E))_\ell$ is a continuous extension of the linear operator $d^{\nabla_{\!0}}$ defined on $\Omega^0_L(M,\End(E))$ and
\begin{align}
\Ad \colon \Gau(E)_{\ell+1} \times \Omega^1_L(M,\End(E))_\ell \rarr \Omega^1_L(M,\End(E))_\ell
\end{align}
is a continuous extension of the left action of the group of gauge transformations $\Gau(E)$ on $\Omega^1_L(M,\End(E))$ defined by \eqref{Ad action} to the appropriate Sobolev completions in the range $\ell > {1 \over 2} \dim M$.

It is easy to see that this action is a smooth mapping of Hilbert manifolds and that, in case $\nabla= \nabla_{\!0} + \alpha \in \eus{A}(E,L)_\ell$ is fixed, the mapping from $\Gau(E)_{\ell+1}$ to $\eus{A}(E,L)_\ell$ given via $\varphi \mapsto \nabla \cdot \varphi$ has the tangent mapping at ${\rm id}_E$ equal to
\begin{align}
  d^\nabla \colon \Omega^0_L(M,\End(E))_{\ell+1} \rarr \Omega^1_L(M,\End(E))_\ell,
\end{align}
where $d^\nabla$ is defined through
\begin{align}
  d^\nabla \!\gamma = d^{\nabla_{\!0}} \gamma + [\alpha,\gamma] \label{covariant}
\end{align}
and $[\cdot\,,\cdot] \colon \Omega^1_L(M,\End(E))_\ell \times \Omega^0_L(M,\End(E))_{\ell+1} \rarr \Omega^1_L(M,\End(E))_\ell$ is a continuous extension of \eqref{lie-bracket} by Sobolev multiplication theorem in the range $\ell > {1 \over 2} \dim M$.
\smallskip

Furthermore, the curvature of an $L$-connection $\nabla= \nabla_{\!0} +\alpha \in \eus{A}(E,L)_\ell$ is defined, using \eqref{connection-curvature2}, by the formula
\begin{align}
  R^\nabla=R^{\nabla_{\!0}+\alpha}=R^{\nabla_{\!0}}+ d^{\nabla_{\!0}}\alpha + {1 \over 2}\,[\alpha,\alpha],
\end{align}
where $\alpha \in \Omega^1_L(M,\End(E))_\ell$, $d^{\nabla_{\!0}} \colon \Omega^1_L(M,\End(E))_\ell \rarr \Omega^2_L(M,\End(E))_{\ell-1}$ is a continuous extension of the linear operator $d^{\nabla_{\!0}}$ defined on $\Omega^1_L(M,\End(E))$ and the bracket $[\cdot\,,\cdot]$ is an extension of \eqref{lie-bracket} to a continuous bilinear mapping $\Omega^1_L(M,\End(E))_\ell \times \Omega^1_L(M,\End(E))_\ell \rarr \Omega^1_L(M,\End(E))_\ell$ in the range $\ell > {1 \over 2} \dim M$.

It is easy to see that $F \colon \eus{A}(E,L)_\ell \rarr \Omega^2_L(M,\End(E))_{\ell-1}$, defined via $F(\nabla)=R^\nabla$, is a smooth mapping of Hilbert manifolds, and the tangent mapping
\begin{gather}
  T_\nabla F \colon \Omega^1_L(M,\End(E))_\ell \rarr \Omega^2_L(M,\End(E))_{\ell-1} \nonumber
  \intertext{is given by}
  (T_\nabla F)(\gamma) =d^{\nabla_{\!0}}\gamma + [\alpha,\gamma] =d^\nabla \!\gamma,
\end{gather}
where $\nabla=\nabla_{\!0} +\alpha \in \eus{A}(E,L)_\ell$ and $\gamma \in \Omega^1_L(M,\End(E))_\ell$.
\medskip

\remark For $\ell > {1 \over 2} \dim M$, we will denote by $\eus{H}(E,L)_\ell$ the set of flat Sobolev $L$-connections. Because $F \colon \eus{A}(E,L)_\ell \rarr \Omega^2_L(M,\End(E))_{\ell-1}$ is a continuous mapping, $\eus{H}(E,L)_\ell$ is a closed subset of $\eus{A}(E,L)_\ell$. Moreover, if we fix some flat $L$-connection $\nabla_{\!0} \in \eus{H}(E,L)$, then
\begin{align}
  \eus{H}(E,L)_\ell=\{\nabla_{\!0} + \alpha;\, \alpha \in \Omega^1_L(M,\End(E))_\ell,\, d^{\nabla_{\!0}}\alpha+ {\textstyle {1\over 2}}[\alpha, \alpha]=0\}.
\end{align}
Furthermore, we need to show that $\eus{H}(E,L)_\ell$ is invariant under the action of the group of gauge transformations $\Gau(E)_{\ell+1}$.
\medskip

\lemma{Let $\nabla=\nabla_{\!0}+\alpha \in \eus{A}(E,L)_\ell$ be a Sobolev $L$-connection on a vector bundle $E$ and let $\varphi \in \Gau(E)_{\ell+1}$ be a gauge transformation, then for $\ell > {1 \over 2} \dim M +1$ we have
\begin{align}
  R^{\nabla^\varphi}=\Ad_{\varphi^{-1}}(R^\nabla),
\end{align}
where $\Ad \colon \Gau(E)_{\ell+1} \times \Omega^2_L(M,\End(E))_{\ell-1} \rarr \Omega^2_L(M,\End(E))_{\ell-1}$ is a continuous extension of \eqref{Ad action} to the appropriate Sobolev spaces using the Sobolev multiplication theorem.}

\proof{If $\nabla=\nabla_{\!0} + \alpha  \in \eus{A}(E,L)_\ell$, then $R^\nabla= R^{\nabla_{\!0}} + d^{\nabla_{\!0}}\alpha+ \alpha \wedge \alpha $. Consider $\varphi \in \Gau(E)_{\ell+1}$, then we have $\nabla^\varphi = \nabla_{\!0} + \varphi^{-1} \wedge (d^{\nabla_{\!0}} \varphi) + \Ad_{\varphi^{-1}}(\alpha)$.
In case we denote $\varphi^{-1} \wedge (d^{\nabla_{\!0}} \varphi) + \Ad_{\varphi^{-1}}(\alpha)$ by $\alpha^\varphi$, then the curvature of the $L$-connection $\nabla^\varphi$ is given by
\begin{align*}
 R^{\nabla^\varphi}= R^{\nabla_{\!0}}+ d^{\nabla_{\!0}}\alpha^\varphi + \alpha^\varphi\! \wedge \alpha^\varphi.
\end{align*}
For $d^{\nabla_{\!0}}\alpha^\varphi$ we can write
\begin{align*}
  d^{\nabla_{\!0}}\alpha^\varphi 
  &=  (d^{\nabla_{\!0}}\varphi^{-1}) \wedge (d^{\nabla_{\!0}} \varphi) + \varphi^{-1} \wedge ((d^{\nabla_{\!0}}\circ d^{\nabla_{\!0}})\,\varphi)\\ & \quad +
(d^{\nabla_{\!0}}\varphi^{-1}) \wedge \alpha \wedge \varphi + \varphi^{-1} \wedge (d^{\nabla_{\!0}} \alpha) \wedge \varphi - \varphi^{-1} \wedge \alpha \wedge (d^{\nabla_{\!0}}\varphi)
\end{align*}
and if we used the fact that $(d^{\nabla_{\!0}} \circ d^{\nabla_{\!0}})\varphi=[R^{\nabla_{\!0}},\varphi]= R^{\nabla_{\!0}} \wedge \varphi - \varphi \wedge R^{\nabla_{\!0}}$, we obtain
\begin{align*}
  d^{\nabla_{\!0}}\alpha^\varphi &= \Ad_{\varphi^{-1}}(R^{\nabla_{\!0}}) + \Ad_{\varphi^{-1}}(d^{\nabla_{\!0}}\alpha) -R^{\nabla_{\!0}} \\ & \quad + (d^{\nabla_{\!0}}\varphi^{-1}) \wedge (d^{\nabla_{\!0}} \varphi) +
(d^{\nabla_{\!0}}\varphi^{-1}) \wedge \alpha \wedge \varphi - \varphi^{-1} \wedge \alpha \wedge (d^{\nabla_{\!0}}\varphi).
\end{align*}
On the other hand, for $\alpha^\varphi \!\wedge \alpha^\varphi$ we get
\begin{align*}
  \alpha^\varphi \!\wedge \alpha^\varphi 
  & = \varphi^{-1} \wedge (d^{\nabla_{\!0}} \varphi) \wedge \varphi^{-1} \wedge (d^{\nabla_{\!0}} \varphi) + \varphi^{-1} \wedge (d^{\nabla_{\!0}} \varphi) \wedge \Ad_{\varphi^{-1}}(\alpha)\\ & \quad + \Ad_{\varphi^{-1}}(\alpha) \wedge \varphi^{-1} \wedge (d^{\nabla_{\!0}} \varphi) + \Ad_{\varphi^{-1}}(\alpha) \wedge \Ad_{\varphi^{-1}}(\alpha)
\end{align*}
and if used the fact that $\varphi \wedge (d^{\nabla_{\!0}} \varphi^{-1})= - (d^{\nabla_{\!0}}\varphi) \wedge \varphi^{-1}$, then we get
\begin{align*}
  \alpha^\varphi \!\wedge \alpha^\varphi
  & = -(d^{\nabla_{\!0}} \varphi^{-1}) \wedge (d^{\nabla_{\!0}} \varphi) - (d^{\nabla_{\!0}} \varphi^{-1}) \wedge \alpha \wedge \varphi + \varphi^{-1} \wedge \alpha \wedge (d^{\nabla_{\!0}} \varphi) + \Ad_{\varphi^{-1}}(\alpha \wedge \alpha).
\end{align*}
Together we obtain
\begin{align*}
  R^{\nabla^\varphi} = \Ad_{\varphi^{-1}}(R^{\nabla_{\!0}}) + \Ad_{\varphi^{-1}}(d^{\nabla_{\!0}}\alpha) + \Ad_{\varphi^{-1}}(\alpha \wedge \alpha)= \Ad_{\varphi^{-1}}(R^\nabla)
\end{align*}
and this finishes the proof.}

Similarly to the smooth case, we define the notion of the irreducibility of Sobolev $L$-connections. For any Sobolev $L$-connection $\nabla$ on a vector bundle $E$ the stabilizer subgroup $\Gau(E)_{\ell+1}^\nabla$ of $\nabla$ contains the subgroup $\Gau(E)^0$. In case $\Gau(E)_{\ell+1}^\nabla = \Gau(E)^0$, then the Sobolev $L$-connection $\nabla$ is called \emph{irreducible}; otherwise, $\nabla$ is called \emph{reducible}. We can prove the following characterization of irreducibility.
\medskip

\lemma{\label{irreducibility - sobolev}Let $\nabla= \nabla_{\!0} + \alpha \in \eus{A}(E,L)_\ell$ be a Sobolev $L$-connection on a vector bundle $E$. Then for $\ell > {1 \over 2} \dim M$ the following conditions are equivalent:
\begin{enum}
  \item[i)] $\Gau(E)_{\ell+1}^\nabla= \Gau(E)^0$,
  \item[ii)] $\ker d^\nabla = \gau(E)^0$.
\end{enum}}

\proof{The proof goes along the similar line as in Lemma \ref{lemma-irred}. Let $\nabla=\nabla_{\!0} + \alpha$ be an $L$-connection and consider a gauge transformation $\varphi \in \Gau(E)_{\ell+1}$. Note that the condition $\nabla \cdot \varphi = \nabla$ means that $\varphi^{-1} \wedge d^{\nabla_{\!0}} \varphi + \Ad_{\varphi^{-1}}(\alpha)= \alpha$. If we multiply this equation by $\varphi$ from the left, we obtain $d^{\nabla_{\!0}} \varphi +[\alpha,\varphi] =0$ and, using \eqref{covariant}, we have $d^\nabla\! \varphi =0$. Therefore, $\varphi \in \Gau(E)_{\ell+1}^\nabla$ if and only if $d^\nabla \!\varphi =0$ and $\varphi \in \Gau(E)_{\ell+1}$.

Suppose that $\varphi \in \Gau(E)_{\ell+1}^\nabla$, then $d^\nabla\! \varphi =0$ and, provided that the condition ii) holds, we obtain $\varphi = f\cdot{\rm id}_E$ for some $f \in C^\infty(M,\mathbb{K})$ satisfying $d_Lf=0$. Thus, we get $\Gau(E)_{\ell+1}^\nabla \subset \Gau(E)^0$ and because the converse inclusion is trivial, we have proved ii) $\Rightarrow$ i).

Now assume that $\varphi \in \ker d^\nabla$. For $c \in \mathbb{K}$ such that $\|\varphi\|_{\ell+1} < |c|\cdot\|{\rm id}_E\|_{\ell+1}$ we obtain that
$c\cdot {\rm id}_E+ \varphi \in \Gau(E)_{\ell+1}$. Furthermore $d^\nabla\!(c\cdot {\rm id}_E + \varphi)=0$ hence, from the previous consideration, we have $c\cdot {\rm id}_E+ \varphi \in \Gau(E)_{\ell+1}^\nabla$. Moreover, if we suppose $\Gau(E)_{\ell+1}^\nabla = \Gau(E)^0$, we obtain $\ker d^\nabla \subset \gau(E)^0$. Converse inclusion is trivial, so we have proved the converse implication.}

We will denote by $\eus{A}^*\!(E,L)_\ell$ the subset of $\eus{A}(E,L)_\ell$ consisting of all irreducible $L$-connections and similarly by $\eus{H}^*\!(E,L)_\ell$ the subset of $\eus{H}(E,L)_\ell$ containing all irreducible flat $L$-connections for $\ell > {1\over 2} \dim M$. It follows from the fact $\Gau(E)^{\nabla^\varphi}_{\ell+1}=\varphi^{-1} \cdot \Gau(E)^\nabla_{\ell+1} \cdot \varphi$ that the irreducibility of $L$-connection is invariant under gauge transformations.

In analogy with \eqref{moduli-konexe}, \eqref{moduli-flat-konexe}, \eqref{moduli-ir-konexe} and \eqref{moduli-ir-flat-konexe} we define the \emph{moduli space}
\begin{align}
    \eus{B}(E,L)_\ell =  \lower-0.5ex\hbox{$\eus{A}(E,L)_\ell$}\mkern-2mu \big/ \mkern-2mu \lower0.5ex\hbox{$\Gau(E)_{\ell+1}$} \label{sobolev-moduli-konexe} \qquad \text{and} \qquad
    \eus{B}^*\!(E,L)_\ell =  \lower-0.5ex\hbox{$\eus{A}^*\!(E,L)_\ell$}\mkern-2mu \big/ \mkern-2mu \lower0.5ex\hbox{$\Gau(E)_{\ell+1}$}
\end{align}
of $L$-connections and irreducible $L$-connections on a vector bundle $E$ for $\ell > {1\over 2} \dim M$ and similarly the \emph{moduli space}
\begin{align}
    \eus{M}(E,L)_\ell =  \lower-0.5ex\hbox{$\eus{H}(E,L)_\ell$}\mkern-2mu \big/ \mkern-2mu \lower0.5ex\hbox{$\Gau(E)_{\ell+1}$} \label{sobolev-moduli-flat-konexe}
    \qquad \text{and} \qquad
    \eus{M}^*\!(E,L)_\ell =  \lower-0.5ex\hbox{$\eus{H}^*\!(E,L)_\ell$}\mkern-2mu \big/ \mkern-2mu \lower0.5ex\hbox{$\Gau(E)_{\ell+1}$}
\end{align}
of flat $L$-connections and irreducible flat $L$-connections on a vector bundle $E$ for $\ell > {1 \over 2} \dim M +1$. Each of these is assumed to have the quotient topology and in the next we will show that $\eus{B}^*\!(E,L)_\ell$ is open in $\eus{B}(E,L)_\ell$ and that $\eus{M}^*\!(E,L)_\ell$ is open in $\eus{M}(E,L)_\ell$. Furthermore, we will denote by
\begin{align}
p_\ell \colon   \eus{A}(E,L)_\ell \rarr \eus{B}(E,L)_\ell
\end{align}
possibly by
\begin{align}
p^*_\ell \colon   \eus{A}^*\!(E,L)_\ell \rarr \eus{B}^*\!(E,L)_\ell
\end{align}
the canonical projection.
\medskip

For $\alpha \in \Omega^1_L(M,\End(E))$ the zero order operator $\ad(\alpha)^* \colon \Omega^1_L(M,\End(E)) \rarr \Omega^0_L(M,\End(E))$ defined as a formal adjoint of $\ad(\alpha) \colon \Omega^0_L(M,\End(E)) \rarr \Omega^1_L(M,\End(E))$, $\ad(\alpha)(\gamma)=[\alpha,\gamma]$, with respect to the Euclidean (Hermitian) metric on $\End(E)$ given by $(f_1,f_2) \mapsto \tr(f_1 \circ f_2^*)$, yields the mapping
\begin{gather}
 \Omega^1_L(M,\End(E)) \times \Omega^1_L(M,\End(E)) \rarr \Omega^0_L(M,\End(E)), \nonumber \\ (\alpha,\beta) \mapsto \ad(\alpha)^*(\beta),
\end{gather}
which is $C^\infty(M,\mathbb{K})$-antilinear in the first component and $C^\infty(M,\mathbb{K})$-linear in the second component. This mapping can be extended through the Sobolev multiplication theorem to a continuous antilinear-linear mapping
\begin{align*}
  \Omega^1_L(M,\End(E))_\ell \times \Omega^1_L(M,\End(E))_k \rarr \Omega^0_L(M,\End(E))_k
\end{align*}
in the range $\ell \geq k > {1 \over 2} \dim M$. Hence the mapping $\ad(\alpha)^* \colon \Omega^1_L(M,\End(E))_k \rarr \Omega^0_L(M,\End(E))_k$ is continuous for all $\alpha \in \Omega^1_L(M,\End(E))_\ell$. Similarly, the mapping
\begin{gather}
 \Omega^1_L(M,\End(E)) \times \Omega^0_L(M,\End(E)) \rarr \Omega^1_L(M,\End(E)), \nonumber \\ (\alpha,\beta) \mapsto \ad(\alpha)(\beta)
\end{gather}
can be extended by the Sobolev multiplication theorem to a continuous bilinear mapping
\begin{align*}
  \Omega^1_L(M,\End(E))_\ell \times \Omega^0_L(M,\End(E))_k \rarr \Omega^1_L(M,\End(E))_k
\end{align*}
in the range $\ell \geq k > {1 \over 2} \dim M$. Therefore, for $\nabla= \nabla_{\!0} + \alpha \in \eus{A}(E,L)_\ell$ we may write
\begin{align}
d^\nabla=d^{\nabla_{\!0}}+ \ad(\alpha)\circ i,
\end{align}
where $i \colon \Omega^0_L(M,\End(E))_{k+1} \rarr \Omega^0(M,\End(E))_k$ is a compact embedding. Furthermore, we define
\begin{align}
\delta^\nabla \colon \Omega^1_L(M,\End(E))_k \rarr \Omega^0_L(M,\End(E))_{k-1}
\end{align}
through
\begin{align}
\delta^\nabla= \delta^{\nabla_{\!0}}+ i \circ \ad(\alpha)^*,
\end{align}
where $i \colon \Omega^0_L(M,\End(E))_k \rarr \Omega^0_L(M,\End(E))_{k-1}$ is a compact embedding and $\delta^{\nabla_{\!0}}$ is a continuous extension of a formal adjoint of $d^{\nabla_{\!0}}$ with respect to the Euclidean (Hermitian) metric on $\End(E)$ in the range $k > {1 \over 2} \dim M$.
\medskip

\lemma{\label{lemma-injectivity}For $\ell >{1 \over 2} \dim M$ the natural mapping
\begin{align}
j_\ell \colon \eus{B}(E,L) \rarr \eus{B}(E,L)_\ell
\end{align}
is injective.}

\proof{Let $\nabla=\nabla_{\!0} + \alpha$ and $\nabla'=\nabla_{\!0} + \alpha'$ be smooth $L$-connections, and suppose we have a gauge transformation $\varphi \in \Gau(E)_{\ell+1}$ satisfying $\nabla'=\nabla \cdot \varphi$, then for the injectivity of $j_\ell$ it suffices to show that $\varphi$ is smooth. If we denote $\beta=\alpha' - \alpha$, then the requirement $\nabla'=\nabla \cdot \varphi$ is equivalent to $d^\nabla \! \varphi = \varphi \wedge \beta$ and we have
\begin{align*}
 \Delta(\varphi) = (\delta^\nabla \!\circ d^\nabla)(\varphi)=\delta^\nabla\!(\varphi \wedge \beta).
\end{align*}
In case $k \geq \ell+1$, then $\varphi \in \Omega^0_L(M,\End(E))_k$ implies, by the Sobolev multiplication theorem, that $\varphi \wedge \beta \in \Omega^1_L(M,\End(E))_k$, because $\beta$ is smooth. Since $\nabla$ is a smooth $L$-connection, the term on the right hand side in the equation above belongs to $\Omega^0_L(M,\End(E))_{k-1}$, and the elliptic regularity theorem, applied to the elliptic operator $\Delta$, gives $\varphi \in \Omega^0_L(M,\End(E))_{k+1}$. Using the induction on $k$, we get $\varphi \in \Omega^0_L(M,\End(E))_k$ for all $k \geq \ell+1$. From the Rellich's lemma it follows that $\varphi$ is smooth.}

\lemma{\label{lemma-fredholm} Let $\nabla=\nabla_{\!0}+\alpha \in \eus{A}(E,L)_\ell$ be an $L$-connection, then the operator
\begin{align}
  \delta^\nabla \!\circ d^\nabla \colon \Omega^0_L(M,\End(E))_{k+1} \rarr \Omega^0_L(M,\End(E))_{k-1}
\end{align}
is a Fredholm operator for $\ell \geq k > {1 \over 2} \dim M$ and $\ind (\delta^\nabla \!\circ d^\nabla) =0$.}

\proof{For $\nabla=\nabla_{\!0}+\alpha$, we may write $\delta^\nabla \!\circ d^\nabla = (\delta^{\nabla_{\!0}}+ i \circ \ad(\alpha)^*) \circ (d^{\nabla_{\!0}} + \ad(\alpha)\circ i)$. Because $\ad(\alpha)\circ i$ and $i \circ \ad(\alpha)^*$ are compact operators,
\begin{align*}
i \circ \ad(\alpha)^*  \circ d^{\nabla_{\!0}} + \delta^{\nabla_{\!0}} \circ \ad(\alpha) \circ i + i \circ \ad(\alpha)^* \circ \ad(\alpha) \circ i
\end{align*}
is also a compact operator.

The rest of the proof is to verify that $\delta^{\nabla_{\!0}} \circ d^{\nabla_{\!0}}$ is a Fredholm operator. It is enough to show that $\delta^{\nabla_{\!0}} \circ d^{\nabla_{\!0}} \colon \Omega^0_L(M,\End(E)) \rarr \Omega^0_L(M,\End(E))$ is an elliptic operator, i.e., that the principal symbol $\sigma_2(\delta^{\nabla_{\!0}} \circ d^{\nabla_{\!0}})(\xi_x) \colon \End(E)_x \rarr \End(E)_x$ is an isomorphism of vector spaces for all $x\in M$ and $0 \neq \xi_x \in T^*_x\!M$. Obviously, we may write
\begin{align*}
\sigma_2(\delta^{\nabla_{\!0}} \circ d^{\nabla_{\!0}})(\xi_x)= \sigma_1(\delta^{\nabla_{\!0}})(\xi_x) \circ \sigma_1(d^{\nabla_{\!0}})(\xi_x)=-(\sigma_1(d^{\nabla_{\!0}})(\xi_x))^*\! \circ \sigma_1(d^{\nabla_{\!0}})(\xi_x)
\end{align*}
and this is an isomorphism if and only if $\sigma_1(d^{\nabla_{\!0}})(\xi_x)$ is injective. But $\sigma_1(d^{\nabla_{\!0}})(\xi_x)=a^*\!(\xi_x)\, \otimes $, i.e., the symbol is the tensor multiplication by $a^*\!(\xi_x)$, therefore it is injective if and only if $a^*\!(\xi_x) \neq 0$. For that reason, $\sigma_2(\delta^{\nabla_{\!0}}\circ d^{\nabla_{\!0}})(\xi_x)$ is an isomorphism for all $x \in M$ and $0 \neq \xi_x \in T^*_x\!M$ if and only if $a^*\!(\xi_x) \neq 0$. However, the last condition is equivalent to the ellipticity condition for the Lie algebroid $L$. Therefore, the operator $\delta^{\nabla_{\!0}} \circ d^{\nabla_{\!0}}$ is a Fredholm operator. For the index we get
\begin{align*}
  \ind(\delta^\nabla \! \circ d^\nabla) = \ind(\delta^{\nabla_{\!0}} \circ d^{\nabla_{\!0}})=0,
\end{align*}
where the last equality follows from the fact that the smooth operator $\delta^{\nabla_{\!0}} \circ d^{\nabla_{\!0}}$ is formally selfadjoint.}

\lemma{\label{orhogonal} Let $\nabla=\nabla_{\!0}+ \alpha \in \eus{A}(E,L)_\ell$ be an $L$-connection, then we have an $L^2$-orthogonal decomposition
\begin{align}
  \Omega^1_L(M,\End(E))_\ell = \im d^\nabla \oplus \ker \delta^\nabla
\end{align}
into closed subspaces for $\ell > {1 \over 2} \dim M + 2$.}

\proof{Let $\nabla= \nabla_{\!0} + \alpha \in \eus{A}(E,L)_\ell$ be an $L$-connection. Further, denote $\Delta_{\alpha,k+1} =\delta_{\alpha,k} \circ d_{\alpha,k+1}$ the corresponding Fredholm operator, as follows from Lemma \ref{lemma-fredholm}, for $\ell \geq k > {1 \over 2}\dim M$.
Therefore, $\dim \ker \Delta_{\alpha,k+1} < \infty$, $\dim \coker \Delta_{\alpha,k+1} < \infty$  and $\im \Delta_{\alpha,k+1}$ is a closed subspace.
Thus, we get an $L^2$-orthogonal (not $L^2_{k+1}$) decomposition $\Omega^0_L(M,\End(E))_{k+1} = \ker \Delta_{\alpha,k+1} \oplus (\ker \Delta_{\alpha,k+1})^\bot$ into closed subspaces of $\Omega^0_L(M,\End(E))_{k+1}$.

For $\beta \in \Omega^0_L(M,\End(E))_{\ell-1}$ and $\gamma \in \Omega^1_L(M,\End(E))_\ell$ we have $(d_{\alpha,\ell-1}\beta,\gamma)=(\beta, \delta_{\alpha,\ell}\gamma)$, hence we obtain that $\im \delta_{\alpha,\ell} \subset (\ker d_{\alpha,\ell-1})^\bot$. Further, because $\im \Delta_{\alpha,\ell+1} \subset \im \delta_{\alpha,\ell} \subset (\ker d_{\alpha,\ell-1})^\bot = (\ker \Delta_{\alpha,\ell-1})^\bot$ and $\ker \Delta_{\alpha,\ell+1} \subset \ker \Delta_{\alpha,\ell-1}$, we get
\begin{align*}
\infty > \dim \coker \Delta_{\alpha,\ell+1} \geq \dim \ker \Delta_{\alpha,\ell-1} \geq \dim \ker \Delta_{\alpha,\ell+1}.
\end{align*}
But $\ind \Delta_{\alpha,\ell+1} =0$, hence we obtain $\im \Delta_{\alpha,\ell+1} = \im \delta_{\alpha,\ell}$. Furthermore, $\im \Delta_{\alpha,\ell+1} $ is a closed subspace, thus ${\Delta_{\alpha,\ell+1}}$ from $(\ker \Delta_{\alpha,\ell+1})^\bot$ to $\im \Delta_{\alpha,\ell+1}$ is a bijective  continuous linear operator between Banach spaces. Therefore, using the Banach's open mapping theorem, there exists a continuous linear operator $G_{\alpha,\ell+1} \colon \im \Delta_{\alpha,\ell+1} \rarr (\ker \Delta_{\alpha,\ell+1})^\bot$ satisfying $\Delta_{\alpha,\ell+1} \circ G_{\alpha,\ell+1}= {\rm id}_{|\im \Delta_{\alpha,\ell+1}}$. Moreover, $d_{\alpha,\ell+1} \circ G_{\alpha,\ell+1} \circ \delta_{\alpha,\ell}$ in a continuous linear operator.
Finally, if $\beta \in \Omega^1_L(M,\End(E))_\ell$, then we have
\begin{align*}
  \beta = (d_{\alpha,\ell+1} \circ G_{\alpha,\ell+1} \circ \delta_{\alpha,\ell})\beta + (\beta - (d_{\alpha,\ell+1} \circ G_{\alpha,\ell+1} \circ \delta_{\alpha,\ell}) \beta)
\end{align*}
with $\delta_{\alpha,\ell}(\beta - (d_{\alpha,\ell+1} \circ G_{\alpha,\ell+1} \circ \delta_{\alpha,\ell}) \beta)=0$. So we are done.}

\lemma{The set of all irreducible Sobolev $L$-connections $\eus{A}^*\!(E,L)_\ell$ is an open subset of $\eus{A}(E,L)_\ell$ for $\ell > {1 \over 2} \dim M$.}

\proof{Let $\nabla = \nabla_{\!0} + \alpha$ be an $L$-connection.
From Lemma \ref{lemma-fredholm} it follows that
\begin{align*}
\Delta_\alpha= \delta^\nabla\! \circ d^\nabla \colon \Omega^0_L(M,\End(E))_{\ell+1} \rarr \Omega^0_L(M,\End(E))_{\ell-1}
\end{align*}
 is a Fredholm operator. Moreover, the mapping
\begin{align*}
  \eus{A}(E,L)_\ell \rarr \mcal{L}(\Omega^0_L(M,\End(E))_{\ell+1}, \Omega^0_L(M,\End(E))_{\ell-1})
\end{align*}
given by $\nabla_{\!0} + \alpha \mapsto \Delta_\alpha$ is a continuous family of Fredholm operators, hence
\begin{align*}
  \nabla_{\!0} + \alpha \mapsto \dim \ker \Delta_\alpha
\end{align*}
is an upper semi-continuous real function on $\eus{A}(E,L)_\ell$, see \cite{HormanderIII1985}. Because we have $\ker d^\nabla = \ker \Delta_\alpha$ and $\dim \ker d^\nabla \geq \dim\gau(E)^0$, the upper semi-continuity implies that $\eus{A}^*\!(E,L)_\ell$ is an open subset of $\eus{A}(E,L)_\ell$.}

\remark We have just proved that $\eus{A}^*\!(E,L)_\ell$ is an open subset of $\eus{A}(E,L)_\ell$. Because $\eus{B}(E,L)_\ell$ is assumed to have the quotient topology and $p^{-1}_\ell(\eus{B}^*\!(E,L)_\ell)= \eus{A}^*\!(E,L)_\ell$, we get that  $\eus{B}^*\!(E,L)_\ell$ is an open subset of $\eus{B}(E,L)_\ell$.
\medskip

Now, for $\nabla=\nabla_{\!0} +\alpha  \in \eus{A}(E,L)_\ell$ and $\veps > 0$ we consider the Hilbert submanifold
\begin{align}
  \mcal{O}_{\alpha,\veps}=\{\nabla_{\!0}+\alpha+ \beta;\, \beta \in \Omega^1_L(M,\End(E))_\ell,\, \delta^\nabla\!\beta=0,\, \|\beta\|_\ell < \veps\}
\end{align}
of the Hilbert manifold $\eus{A}(E,L)_\ell$. Because $\mcal{O}_{\alpha,\veps}$ is a Hilbert manifold modeled on $\ker \delta^\nabla$, we obtain
\begin{align}
  T_\nabla(\mcal{O}_{\alpha,\veps})= \ker \delta^\nabla.
\end{align}
First note that if $\nabla \in \eus{A}^*\!(E,L)_\ell$, then we may take $\veps$ small enough to ensure $\mcal{O}_{\alpha,\veps} \subset \eus{A}^*\!(E,L)_\ell$, since $\eus{A}^*\!(E,L)_\ell$ is open in $\eus{A}(E,L)_\ell$. Next, we define the \emph{reduced group of gauge transformations} $\Gau(E)_{\ell+1}^{\rm r}$ by
\begin{align}
  \Gau(E)_{\ell+1}^{\rm r}= \lower-0.5ex\hbox{$\Gau(E)_{\ell+1}$}\mkern-2mu \big/ \mkern-2mu \lower0.5ex\hbox{$\Gau(E)^0$}. \label{reduced-gauge-group-sobolev}
\end{align}
Because $\Gau(E)^0$ is a closed normal Hilbert--Lie subgroup of $\Gau(E)_{\ell+1}$, Theorem \ref{normal-subgroup} bellow implies that the reduced group of gauge transformations is a Hilbert--Lie group with the Lie algebra
\begin{align}
  \gau(E)_{\ell+1}^{\rm r}= \lower-0.5ex\hbox{$\gau(E)_{\ell+1}$}\mkern-2mu \big/ \mkern-2mu \lower0.5ex\hbox{$\gau(E)^0$},
\end{align}
where the Lie bracket descents from the one on $\gau(E)_{\ell+1}$. Moreover, if
\begin{align}
  q \colon \Gau(E)_{\ell+1} \rarr \Gau(E)_{\ell+1}^{\rm r}=\lower-0.5ex\hbox{$\Gau(E)_{\ell+1}$}\mkern-2mu \big/ \mkern-2mu \lower0.5ex\hbox{$\Gau(E)^0$}
\end{align}
denotes the canonical projection, then $q$ is a smooth mapping and any mapping $f \colon \Gau(E)_{\ell+1}^{\rm r} \rarr X$, where $X$ is a smooth Banach manifold, is smooth if and only if $f \circ q \colon \Gau(E)_{\ell+1} \rarr X$ is smooth.
\medskip

\theorem{\label{normal-subgroup} Let $G$ be a Banach--Lie group over $\mathbb{K}$ with the Lie algebra $\mfrak{g}$ and suppose that $N$ is a closed normal Banach--Lie subgroup over $\mathbb{K}$ of $G$ with the Lie algebra $\mfrak{n}$. Then $G/N$ is a Banach--Lie group over $\mathbb{K}$ with the Lie algebra $\mfrak{g}/\mfrak{n}$ in a unique way such that the quotient mapping $q \colon G \rarr G/N$ is smooth. Moreover, for any Banach manifold $X$ a mapping $f \colon G/N \rarr X$ is smooth if and only if $f \circ q$ is smooth.}

\proof{See \cite{Glockner_Neeb2003}, \cite{Glockner2002} and \cite{Dupre_Glazebrook2000}.}

\theorem{\label{moduli-ir-hilbert}$\eus{B}^*\!(E,L)_\ell$ is a locally Hausdorff Hilbert manifold and $p^*_\ell \colon \eus{A}^*\!(E,L)_\ell \rarr \eus{B}^*\!(E,L)_\ell$ is a principal $\Gau(E)_{\ell+1}^{\rm r}$-bundle for $\ell > {1 \over 2} \dim M+2$.}

\proof{Let $\nabla=\nabla_{\!0} + \alpha$ be an irreducible $L$-connection. Consider the smooth mapping of Hilbert manifolds
\begin{gather*}
  \Psi_\nabla \colon \mcal{O}_{\alpha,\veps} \times \Gau(E)_{\ell+1}^{\rm r} \rarr \eus{A}^*\!(E,L)_\ell, \nonumber \\
  \Psi_\nabla(\nabla_{\!0}+\alpha+\beta,[\varphi])= (\nabla_{\!0}+\alpha+\beta) \cdot \varphi,
\end{gather*}
then the tangent mapping at $(\nabla,[{\rm id}_E])$ equals to
\begin{gather*}
  T_{(\nabla,[{\rm id}_E])} \Psi_\nabla \colon \ker \delta^\nabla \oplus \gau(E)^{\rm r}_{\ell+1} \rarr \Omega^1_L(M,\End(E))_\ell, \nonumber \\
  (T_{(\nabla,[{\rm id}_E])} \Psi_\nabla) (\beta, [\gamma]) = d^\nabla \!\gamma + \beta.
\end{gather*}
From Lemma~\ref{orhogonal} it follows that $T_{(\nabla,[{\rm id}_E])} \Psi_\nabla$ is surjective. Moreover, because $\nabla$ is assumed to be an irreducible $L$-connection, we obtain, using  Lemma~\ref{irreducibility - sobolev}, that $T_{(\nabla,[{\rm id}_E])}\Psi_\nabla$ is injective. Hence by the Banach's open mapping theorem $T_{(\nabla,[{\rm id}_E])} \Psi_\nabla$ is an isomorphism. Therefore, the inverse function theorem for Banach manifolds implies that $\Psi_\nabla$ is a local diffeomorphism near $(\nabla,[{\rm id}_E])$. Consequently, there is an open neighborhood $\mcal{V}_\alpha$ of $\nabla$ in $\eus{A}^*\!(E,L)_\ell$ and an open neighborhood $\mcal{N}_{[{\rm id}_E]}$ of $[{\rm id}_E]$ in $\Gau(E)_{\ell+1}^{\rm r}$ such that
\begin{align*}
  \Psi_\nabla \colon \mcal{O}_{\alpha,\veps} \times \mcal{N}_{[{\rm id}_E]} \rarr \mcal{V}_\alpha
\end{align*}
is a diffeomorphism for sufficiently small $\veps > 0$.

Next we will show that, for $\veps$ small enough, the mapping $p_{\alpha,\veps}={p^*_\ell}_{|\mcal{O}_{\alpha,\veps}} \colon \mcal{O}_{\alpha,\veps} \rarr \eus{B}^*\!(E,L)_\ell$ is injective. We have to show that if for two elements $\nabla_{\!0}+ \alpha + \beta_1, \nabla_{\!0}+\alpha+\beta_2 \in \mcal{O}_{\alpha,\veps}$ there exists a gauge transformation $\varphi \in \Gau(E)_{\ell+1}$ satisfying
\begin{align*}
  (\nabla_{\!0}+\alpha+\beta_1) \cdot \varphi= \nabla_{\!0} + \alpha +\beta_2,
\end{align*}
then $\beta_1= \beta_2$. First observe that the previous equation is equivalent to
\begin{align*}
  d^\nabla \!\varphi = \varphi \wedge \beta_2 - \beta_1 \wedge \varphi.
\end{align*}
Further, because $\Omega^0_L(M,\End(E))_{\ell+1}= \ker d^\nabla \oplus (\ker d^\nabla)^\bot$ is an $L^2$-orthogonal decomposition into closed subspaces, we can write $\varphi= f \cdot {\rm id}_E + \varphi_0$, where $f\in C^\infty(M,\mathbb{K})$, $d_Lf=0$ and $\varphi_0 \in (\ker d^\nabla)^\bot$. Moreover, $\im d^\nabla$ is a closed subspace of $\Omega^1_L(M,\End(E))_\ell$, hence we obtain by the Banach's open mapping theorem that
\begin{align*}
  d^\nabla \colon (\ker d^\nabla)^\bot \rarr \im d^\nabla
\end{align*}
is an isomorphism of Hilbert spaces. Therefore, it is a lower bounded operator, i.e., there exists a positive constant $c_1$ such that
\begin{align*}
   \|d^\nabla \!\psi\|_\ell \geq c_1 \|\psi\|_{\ell+1}
\end{align*}
for all $\psi \in (\ker d^\nabla)^\bot$. Denote by $M_i$, for $i=1,2,\dots,n$, connected components of $M$. If we consider $\psi_1, \psi_2 \in \Omega^0_L(M,\End(E))_{\ell+1}$ such that $\esssupp \psi_1 \subset M_{i_1}$ and $\esssupp \psi_2 \subset M_{i_2}$, then it is easy to see that $(\psi_1,\psi_2)=0$ provided $i_1 \neq i_2$. This gives an $L^2_{\ell+1}$-orthogonal decomposition of $\Omega^0_L(M,\End(E))_{\ell+1}$ into closed subspaces $\Omega^0_L(M,\End(E))^i_{\ell+1}$ defined by
\begin{align*}
  \Omega^0_L(M,\End(E))^i_{\ell+1}=\{\psi \in \Omega^0_L(M,\End(E))_{\ell+1};\, \esssupp \psi \subset M_i\}
\end{align*}
for $i=1,2,\dots,n$. For $\psi \in \Omega^0_L(M,\End(E))_{\ell+1}$ the corresponding elements of $\Omega^0_L(M,\End(E))^i_{\ell+1}$ will be denoted by $\psi^i$.

Therefore, we may write
\begin{align*}
  c_1 \|\varphi^i_0\|_{\ell+1} \leq \|d^\nabla \!\varphi^i_0 \|_\ell= \|d^\nabla \! \varphi^i\|_\ell= \|\varphi^i \wedge \beta_2 - \beta_1 \wedge \varphi^i\|_\ell \leq 2\veps c_0 \, (\|f\cdot{\rm id}^i_E\|_{\ell+1} + \|\varphi^i_0\|_{\ell+1}),
\end{align*}
where we used the fact that $\|\psi \wedge \gamma\|_\ell \leq c_0  \|\psi\|_{\ell+1} \|\gamma\|_\ell$ and $\|\gamma \wedge \psi\|_\ell \leq c_0 \|\psi\|_{\ell+1}\|\gamma\|_\ell $ for all $\psi \in \Omega^0_L(M,\End(E))_{\ell+1}$ and $\gamma \in \Omega^1_L(M,\End(E))_\ell$. As a consequence we have
\begin{align*}
  \|\varphi^i_0\|_{\ell+1} \leq {2\veps c_0  \over c_1 - 2 \veps c_0 } \, \|f\cdot{\rm id}^i_E\|_{\ell+1} = {2\veps c_0  \over c_1 - 2 \veps c_0 } \, \|c^i\cdot{\rm id}^i_E\|_{\ell+1}
\end{align*}
for $\veps$ small enough. The last equality follows from the fact that $f=c^i$ for some $c^i \in \mathbb{K}$ on $M^i$ because $f$ is a locally constant function on $M$ as we have from Lemma \ref{lemma-kernel}. Provided that $c^i=0$, we obtain immediately $\|\varphi^i_0\|_{\ell+1}=0$, therefore $\varphi^i=0$ and this is a contradiction. Furthermore, because $f \in \Gau(E)^0$, we get
\begin{align*}
  \|f^{-1}\cdot \varphi - {\rm id}_E\|_{\ell+1} & = \|f^{-1}\cdot\varphi_0\|_{\ell+1}  \leq  {2\veps c_0 \over c_1 - 2 \veps c_0} \, \|f^{-1}f \cdot {\rm id}_E\|_{\ell+1} = {2\veps c_0 \over c_1 - 2 \veps c_0} \, \|{\rm id}_E\|_{\ell+1}.
\end{align*}
Since $q^{-1}(\mcal{N}_{[{\rm id}_E]})$ is an open set in $\Gau(E)_{\ell+1}$ and ${\rm id}_E \in q^{-1}(\mcal{N}_{[{\rm id}_E]})$, we obtain that for $\veps$ small enough $[\varphi]$ is near $[{\rm id}_E]$ in $\Gau(E)_{\ell+1}^{\rm r}$, i.e., $[\varphi] \in \mcal{N}_{[{\rm id}_E]}$.
And if we use that $\Psi_\nabla$ is injective, we get $\beta_1 =\beta_2$.
\smallskip

Let $\mcal{U}_\alpha= p^*_\ell(\mcal{O}_{\alpha,\veps})$, then  we have $(p^*_\ell)^{-1}(\mcal{U}_\alpha)=\tilde{R}(\mcal{O}_{\alpha,\veps} \times \Gau(E)_{\ell+1})$, where $\tilde{R} \colon \eus{A}^*\!(E,L)_\ell \times \Gau(E)_{\ell+1} \rarr \eus{A}^*\!(E,L)_\ell$ is the corresponding right action. From the previous considerations, it follows that $(p^*_\ell)^{-1}(\mcal{U}_\alpha)$ is open in  $\eus{A}^*\!(E,L)_\ell$, therefore $\mcal{U}_\alpha$ is open in $\eus{B}^*\!(E,L)_\ell$. Moreover, $p_{\alpha,\veps} \colon \mcal{O}_{\alpha,\veps} \rarr \mcal{U}_\alpha$ is a homeomorphism. The mapping
\begin{gather*}
  \Psi_\nabla \colon \mcal{O}_{\alpha,\veps} \times \Gau(E)_{\ell+1}^{\rm r} \rarr (p^*_\ell)^{-1}(\mcal{U}_\alpha), \nonumber \\
  \Psi_\nabla(\nabla_{\!0}+\alpha+\beta,[\varphi])= (\nabla_{\!0}+ \alpha+ \beta) \cdot \varphi
\end{gather*}
is surjective because $(p^*_\ell)^{-1}(\mcal{U}_{\alpha,\veps})= \tilde{R}(\mcal{O}_{\alpha,\veps} \times \Gau(E)_{\ell+1})$, the injectivity follows from the previous consideration and from the fact that the action of $\Gau(E)_{\ell+1}^{\rm r}$ on $\eus{A}^*\!(E,L)_\ell$ is free. In fact, we will show that $\Psi_\nabla$ is a diffeomorphism of Hilbert manifolds.

For an arbitrary $[\varphi] \in \Gau(E)_{\ell+1}^{\rm r}$ we find an open neighborhood $\mcal{W}_{[\varphi]}$ of $[\varphi]$ such that the mapping ${\Psi_\nabla}_{|\mcal{O}_{\alpha,\veps} \times R_{[\varphi]^{-1}}(\mcal{W}_{[\varphi]})}$
is a diffeomorphism, where $R_{[\varphi]^{-1}}$ is the right translation by $[\varphi]^{-1}$ in $\Gau(E)_{\ell+1}^{\rm r}$. In particular, we can take $\mcal{W}_{[\varphi]}=R_{[\varphi]}(\mcal{N}_{[{\rm id}_E]})$. Therefore, we have
\begin{align*}
  {\Psi_\nabla}_{|\mcal{O}_{\alpha,\veps} \times \mcal{W}_{[\varphi]}}= \tilde{R}_\varphi \circ {\Psi_\nabla} \circ ({\rm id}_{\eus{A}^*\!(E,L)_\ell} \times R_{[\varphi]^{-1}})_{|\mcal{O}_{\alpha,\veps} \times \mcal{W}_{[\varphi]}},
\end{align*}
which is a diffeomorphism.
\medskip

Now, to show that $p^*_\ell \colon \eus{A}^*\!(E,L)_\ell \rarr \eus{B}^*\!(E,L)_\ell$ is a principal $\Gau(E)_{\ell+1}^{\rm r}$-bundle over a Hilbert manifold, we only need to glue together the local charts $u_\alpha \colon \mcal{U}_\alpha \rarr \mcal{O}_{\alpha,\veps}$, $u_\alpha= p^{-1}_{\alpha,\veps}$. Consider the smooth mapping
\begin{align*}
  \psi_\alpha ={\rm pr} \circ \Psi_\nabla^{-1} \colon (p^*_\ell)^{-1}(\mcal{U}_\alpha) \rarr \Gau(E)_{\ell+1}^{\rm r},
\end{align*}
where ${\rm pr} \colon \mcal{O}_{\alpha,\veps} \times \Gau(E)_{\ell+1}^{\rm r} \rarr \Gau(E)_{\ell+1}^{\rm r}$ is the projection. Further, for any $\nabla'= \nabla_{\!0} + \alpha' \in \eus{A}^*\!(E,L)_\ell$ with $p^*_\ell(\nabla_{\!0}+\alpha') \in \mcal{U}_\alpha$ we have
\begin{align*}
  u_\alpha(p^*_\ell(\nabla_{\!0}+\alpha'))= (\nabla_{\!0}+\alpha') \cdot (\psi_\alpha(\nabla_{\!0}+\alpha'))^{-1}.
\end{align*}
Hence it is easy to see that over $u_{\alpha_2}(\mcal{U}_{\alpha_2} \cap \mcal{U}_{\alpha_1})$ we have
\begin{align*}
  (u_{\alpha_1} \circ u_{\alpha_2}^{-1})(\nabla_{\!0}+\alpha_2+\beta)= u_{\alpha_1}(p^*_\ell(\nabla_{\!0}+\alpha_2+\beta))=(\nabla_{\!0}+\alpha_2+\beta) \cdot (\psi_\alpha(\nabla_{\!0}+\alpha_2+\beta))^{-1},
\end{align*}
and this is clearly smooth in $\beta$.}

\section{Moduli spaces -- local model}

In this section, we give a local description of the moduli space $\eus{M}^*\!(E,L)$ of isomorphism classes of irreducible flat $L$-connections around a given point. We will adopt to this situation the Kuranishi's argument used for describing the moduli space of complex structures on a compact manifold and the moduli space of anti-self-dual connections on a compact 4-manifold given by Atiyah, Hitchin and Singer, see \cite{Atiyah_Hitchin_Singer1978}.

The \emph{Kuranishi's description} provides local models of the moduli space, i.e., it gives an explicit description of the germ of the moduli space in a given point. This makes it possible to estimate the dimension of the moduli space and provides a simple smoothness criterium.
\medskip

Let $(L \rarr M, [\cdot\,,\cdot],a)$ be a real (complex) Lie algebroid satisfying the \emph{ellipticity condition} and let $E \rarr M$ be a real (complex) vector bundle. Further, assume that $M$ is a compact manifold. Then to any flat $L$-connection $\nabla$ on $E$ is associated the \emph{fundamental elliptic complex} $\mcal{E}(\nabla)$ playing a cental role in the subsequent discussion.

Consider the sequence of linear differential operators of first order
\begin{align}
\label{deformation-complex}
0 \xrightarrow{\phantom{d^\nabla}} \Omega^0_L(M,\End(E)) \xrightarrow{d^\nabla}
\Omega^1_L(M,\End(E)) \xrightarrow{d^\nabla} \hdots \xrightarrow{d^\nabla} \Omega^r_L(M,\End(E)) \xrightarrow{\phantom{d^\nabla}} 0,
\end{align}
where $r =\rk L$ and $d^\nabla$ is the covariant exterior derivative for the induced $L$-connection $\nabla^{\End(E)}$ on $\End(E)$. Because $R^\nabla=0$ and
\begin{align}
  R^{\nabla^{\End(E)}}\!(\xi,\eta)\,\gamma=[R^\nabla\!(\xi,\eta),\gamma] =[R^\nabla\!,\gamma](\xi,\eta),
\end{align}
where $\xi,\eta \in \mfrak{X}_L(M)$ and $\gamma \in \Omega^0_L(M,\End(E))$, we obtain  $R^{\nabla^{\End(E)}}=0$. Further, using Lemma~\ref{lemma-curvature} and the fact that the Lie algebroid satisfies the ellipticity condition, we get that the sequence \eqref{deformation-complex} of differential operators is an elliptic complex, called the \emph{deformation complex}. We will denote the cohomology of this elliptic compex by $H^i(E,\nabla)$ for $i=0,1,\dots,r$.

Endow $E$ and $L$ with an Euclidean (Hermitian) metric $h_E$ and $h_L$, respectively. This gives an Euclidean (Hermitian) metric on each vector bundle $\Lambda^kL^*\otimes \End(E)$. Furthermore, let $g$ be a Riemannian metric on $M$. Then we have the formal self-adjoint elliptic operators of second order
\begin{align}
  \Delta_i = \delta^\nabla_i \!\circ d^\nabla_i + d^\nabla_{i-1} \circ \delta^\nabla_{i-1} \colon \Omega^i_L(M,\End(E))  \rarr \Omega^i_L(M,\End(E)),
\end{align}
where $\delta^\nabla_i$ is a formal adjoint of $d^\nabla_i$ and $d^\nabla_{-1}, d^\nabla_r$ are zero operators. Besides, the kernel of $\Delta_i$
\begin{align}
  \mcal{H}^i(E,\nabla)=\{\alpha \in \Omega^i_L(M,\End(E));\, \Delta_i \alpha=0\}= \ker d^\nabla_i \cap \ker \delta^\nabla_{i-1}
\end{align}
is a finite dimensional vector space for $i=0,1,\dots,r$ and moreover there exists a natural isomorphism $\mcal{H}^i(E,\nabla) \simeq H^i(E,\nabla)$. Because all cohomology groups are finite dimensional vector spaces, we may define the \emph{index} of $\mcal{E}(\nabla)$ by
\begin{align}
  \ind \mcal{E}(\nabla) = \sum_{i=0}^r (-1)^i \dim H^i(E,\nabla) = \sum_{i=0}^r (-1)^i \dim \ker \Delta_i.
\end{align}
A fundamental result of the Hodge theory for the elliptic complex \eqref{deformation-complex} is the Hodge decomposition theorem, which states that there is an $L^2$-orthogonal decomposition
\begin{align}
  \Omega^i_L(M,\End(E))= \mcal{H}^i(E,\nabla) \oplus \im d^\nabla_{i-1} \oplus \im \delta^\nabla_i.
\end{align}
Furthermore, there exists a unique linear pseudo-differential operator
\begin{align}
  G_i \colon \Omega^i_L(M,\End(E)) \rarr \Omega^i_L(M,\End(E)),
\end{align}
of order $-2$, called the Green's operator associated to $\Delta_i$, satisfying
\begin{align}
  {\rm id}_{\Omega^i_L(M,\End(E))} = H_i + \Delta_i \circ G_i = H_i + G_i \circ \Delta_i \label{identity-flat}
\end{align}
and the following commutation relations
\begin{align}
   H_i \circ G_i = G_i \circ H_i, \quad d^\nabla_i \circ G_i= G_{i+1} \circ d^\nabla_i, \quad  \delta^\nabla_i \circ G_{i+1}= G_i \circ \delta^\nabla_i, \label{identity-comutator}
\end{align}
where $H_i \colon \Omega^i_L(M,\End(E)) \rarr \mcal{H}^i(E,\nabla)$ for $i=0,1,\dots,r$ are $L^2$-orthogonal projections. Further, all the associated operators $d^\nabla_i$, $\delta^\nabla_i$, $\Delta_i$, $G_i$ can be extended to continuous linear operators between appropriate Sobolev completions, e.g.
\begin{align*}
  d^\nabla_{i,k} \colon \Omega^i_L(M,\End(E))_k \rarr \Omega^{i+1}_L(M,\End(E))_{k-1},  \\
  \delta^\nabla_{i-1,k} \colon \Omega^i_L(M,\End(E))_k \rarr \Omega^{i-1}_L(M,\End(E))_{k-1}, \\
  \Delta_{i,k} \colon \Omega^i_L(M,\End(E))_k \rarr \Omega^i_L(M,\End(E))_{k-2}, \\
  G_{i,k} \colon \Omega^i_L(M,\End(E))_k \rarr \Omega^i_L(M,\End(E))_{k+2},
\end{align*}
and note that
\begin{align}
  \ker \Delta_{i,k} = \ker \Delta_i = \mcal{H}^i(E,\nabla).
\end{align}
for $i=0,1,\dots,r$.
\medskip

\remark Note that $H^0(E,\nabla) = \ker \Delta_0 = \ker d^\nabla_0$. Therefore, if $\nabla$ is an irreducible $L$-connection, then, using Lemma \ref{lemma-irred} and Lemma \ref{lemma-kernel}, we get $\dim H^0(E,\nabla) =b_0(M)$, where $b_0(M)$ is the zero Betti number of $M$, otherwise $\dim H^0(E,\nabla) > b_0(M)$.
\medskip

Recall that if we fix some smooth flat $L$-connection $\nabla_{\!0} \in \eus{H}(E,L)$, then the Sobolev completion is defined through
\begin{align}
  \eus{H}(E,L)_\ell=\{\nabla_{\!0}+\alpha;\, \alpha \in \Omega^1_L(M,\End(E))_\ell,\, d^{\nabla_{\!0}}\alpha + \textstyle{{1 \over 2}}\,[\alpha,\alpha]=0\}
\end{align}
for $\ell > {1\over 2} \dim M$. Furthermore, from the previous we know that the curvature
\begin{align}
  F \colon \eus{A}(E,L)_\ell \rarr \Omega^2_L(M,\End(E))_{\ell-1}
\end{align}
given by $F(\nabla_{\!0}+\alpha)=d^{\nabla_{\!0}}\alpha + {1 \over 2}[\alpha,\alpha]$
is a smooth mapping of Hilbert manifolds for $\ell > {1 \over 2} \dim M$ and
\begin{align}
  \eus{H}(E,L)_\ell = F^{-1}(0).
\end{align}
Consider a smooth irreducible flat $L$-connection $\nabla= \nabla_{\!0} + \alpha \in \eus{H}^*\!(E,L)$. From the proof of Theorem~\ref{moduli-ir-hilbert} we have that for the Hilbert submanifold $\mcal{O}_{\alpha,\veps}$ of $\eus{A}^*\!(E,L)_\ell$, for $\veps > 0$ small enough and $\ell > {1\over 2} \dim M+2$, the mapping $p_{\alpha,\veps} = {p^*_\ell}_{|\mcal{O}_{\alpha,\veps}} \colon \mcal{O}_{\alpha,\veps} \rarr \mcal{U}_\alpha \subset \eus{B}^*\!(E,L)_\ell$, where $\mcal{U}_\alpha=p^*_\ell(\mcal{O}_{\alpha,\veps})$, is a homeomorphism onto an open subset of $\eus{B}^*\!(E,L)_\ell$, i.e., $\mcal{O}_{\alpha,\veps}$ is a slice to the $\Gau(E)_{\ell+1}$-orbits of the action of the group of gauge transformations $\Gau(E)_{\ell+1}$ on the set of irreducible $L$-connections $\eus{A}^*\!(E,L)_\ell$.
Furthermore, consider a closed subset
\begin{align}
  \mcal{S}_{\alpha,\veps}=\{\nabla_{\!0}+\alpha + \beta;\, \beta \in \Omega^1_L(M,\End(E))_\ell,\, \delta^\nabla\!\beta=0,\, d^\nabla\!\beta + \textstyle{{1 \over 2}}[\beta,\beta]=0,\, \|\beta\|_\ell < \veps\}
\end{align}
of $\mcal{O}_{\alpha,\veps}$. As $\mcal{S}_{\alpha,\veps} \subset \eus{H}^*\!(E,L)_\ell$, we obtain that $p_{\alpha,\veps} \colon \mcal{S}_{\alpha,\veps} \rarr \mcal{V}_\alpha = \mcal{U}_\alpha \cap \eus{M}^*\!(E,L)_\ell$ is a homeomorphism onto an open subset of $\eus{M}^*\!(E,L)_\ell$ for $\ell > {1 \over 2} \dim M +2$.
\medskip

Suppose that $\ell > \max\{{1 \over 2} \dim M,1\}$. Now, if we apply the Hodge decomposition \eqref{identity-flat} to the element $d^\nabla_1\!\beta + {1 \over 2} [\beta, \beta]$ for some $\beta \in \Omega^1_L(M,\End(E)_\ell$, we obtain
\begin{align*}
  \begin{split}
    d^\nabla_1\!\beta + {\textstyle {1 \over 2}}[\beta, \beta] &= H_2 (d^\nabla_1\!\beta + {\textstyle {1 \over 2}}[\beta, \beta]) + (\delta^\nabla_2 \circ d^\nabla_2 \circ G_2) (d^\nabla_1\! \beta + {\textstyle {1 \over 2}}[\beta, \beta]) \\ & \quad + (d^\nabla_1 \circ \delta^\nabla_1 \circ  G_2) (d^\nabla_1\! \beta + {\textstyle {1 \over 2}}[\beta, \beta]) \\
    &={\textstyle {1 \over 2}}H_2([\beta, \beta]) + {\textstyle {1 \over 2}} (\delta^\nabla_2 \circ d^\nabla_2 \circ G_2)([\beta, \beta]) \\ & \quad + d^\nabla_1\! ((\delta^\nabla_1 \circ G_2 \circ d^\nabla_1)\beta + {\textstyle {1 \over 2}} (\delta^\nabla_1 \circ  G_2) ([\beta, \beta])),
  \end{split}
\end{align*}
where we used that $G_2 \circ d^\nabla_1 = d^\nabla_1 \circ G_1$. Besides, we have
\begin{align*}
    \delta^\nabla_1 \circ G_2 \circ d^\nabla_1&= \delta^\nabla_1 \circ d^\nabla_1 \circ G_1 = \Delta_1 \circ G_1 - d^\nabla_0 \circ \delta^\nabla_0 \circ G_1\\ &= {\rm id}_{\Omega^1_L(M,\End(E))_\ell} - H_1- d^\nabla_0 \circ \delta^\nabla_0 \circ G_1,
\end{align*}
therefore substituting this into the equation above, we get
\begin{align*}
  d^\nabla_1\!\beta + {\textstyle {1 \over 2}}[\beta, \beta] &= {\textstyle {1 \over 2}} H_2([\beta,\beta]) + {\textstyle {1 \over 2}}(\delta^\nabla_2 \circ d^\nabla_2 \circ G_2)([\beta, \beta]) \\
  &\quad + d^\nabla_1\! (\beta + {\textstyle {1 \over 2}} ( \delta^\nabla_1 \circ G_2)([\beta,\beta])).
\end{align*}
From this $L^2$-orthogonal decomposition we have
\begin{align}\label{decomposition}
  d^\nabla_1\! \beta + {1 \over 2}\, [\beta, \beta] = 0 \Longleftrightarrow
  \begin{cases}
    d^\nabla_1\!\big(\beta + {1 \over 2}(\delta^\nabla_1 \circ G_2)([\beta, \beta])\!\big)=0, \\
    (\delta^\nabla_2 \circ d^\nabla_2 \circ G_2)([\beta, \beta]) = 0, \\
    H_2([\beta, \beta])=0.
  \end{cases}
\end{align}
Furthermore, for the smooth irreducible flat $L$-connection $\nabla$ we define the \emph{Kuranishi mapping}
\begin{align*}
K_\nabla \colon \Omega^1_L(M,\End(E))_\ell \rarr \Omega^1_L(M,\End(E))_\ell
\end{align*}
by the formula
\begin{align}
  K_\nabla(\beta)= \beta + {1 \over 2}\, (\delta^\nabla_1 \circ G_2)([\beta, \beta])
\end{align}
for $\beta \in \Omega^1_L(M,\End(E))_\ell$. It is a smooth mapping of Hilbert manifolds with the tangent mapping $T_\beta K_\nabla \colon \Omega^1_L(M,\End(E))_\ell \rarr \Omega^1_L(M,\End(E))_\ell$ at $\beta$ equal to
\begin{align}
  T_\beta K_\nabla(\gamma)= \gamma + (\delta^\nabla_1 \circ G_2)([\beta, \gamma]),
\end{align}
where $\gamma \in \Omega^1_L(M,\End(E))_\ell$. Since $T_0K_\nabla={\rm id}_{\Omega^1_L(M,\End(E))_\ell}$, using the inverse function theorem for Banach manifolds, we immediately obtain that $K_\nabla$ is a local diffeomorphism at $0$. Further, we define the subset
\begin{align}
  \mcal{S}^\ell_\veps=\{\beta \in \Omega^1_L(M,\End(E))_\ell,\, \delta^\nabla_0\!\beta=0,\, d^\nabla_1\!\beta+{\textstyle {1 \over 2}} [\beta,\beta]=0,\, \|\beta\|_\ell < \veps\}
\end{align}
of $\Omega^1_L(M,\End(E))_\ell$ for $\veps > 0$.
\medskip

\lemma{\label{Kuranishi smooth}Let $\ell > \max\{{1\over 2} \dim M,1\}$, then $K_\nabla(\mcal{S}^\ell_\veps) \subset \mcal{H}^1(E,\nabla)$ and $\mcal{S}^\ell_\veps \subset \Omega^1_L(M,\End(E))$.}

\proof{The first observation is easy. It is enough to show that $d^\nabla_1\!( K_\nabla(\beta))=0$ and $\delta^\nabla_0\!(K_\nabla(\beta))=0$ for $\beta \in \mcal{S}^\ell_\veps$, since $\mcal{H}^1(E,\nabla) = \ker d^\nabla_1 \cap \ker \delta^\nabla_0$. But we have $\delta^\nabla_0\!(K_\nabla (\beta))=\delta^\nabla_0\!\beta=0$ and using \eqref{decomposition} we obtain $d^\nabla_1\!(K_\nabla (\beta))=d^\nabla_1\! (\beta+ {1 \over 2} (\delta^\nabla_1 \circ G_2)([\beta, \beta]))=0$.

Now consider $\beta \in \mcal{S}^\ell_\veps$ and assume that $\beta \in \Omega^1_L(M,\End(E))_k$ for $k \geq \ell$. Because $\Delta_1(K_\nabla(\beta))=0$, we get
\begin{align*}
  \Delta_1\beta=-{1 \over 2}\, (\Delta_1 \circ \delta^\nabla_1 \circ G_2)([\beta,\beta]).
\end{align*}
The term on the right hand side in the equation above belongs to $\Omega^1_L(M,\End(E))_{k-1}$, and the elliptic regularity theorem, applied to the elliptic operator $\Delta_1$, gives $\beta \in \Omega^1_L(M,\End(E))_{k+1}$. Using the induction on $k$, we get $\beta \in \Omega^1_L(M,\End(E))_k$ for all $k \geq \ell$. From the Rellich's lemma it follows that $\beta$ is smooth.}

\lemma{\label{open image}For $\ell > {1 \over 2} \dim M +2$ the mapping $j_\ell \colon \eus{M}^*\!(E,L) \rarr \eus{M}^*\!(E,L)_\ell$ is injective and has an open image.}

\proof{The injectivity of $j_\ell$ follows from Lemma \ref{lemma-injectivity} and the fact that $j_\ell(\eus{M}^*\!(E,L)) \subset \eus{M}^*\!(E,L)_\ell$.
Further, let $\nabla=\nabla_{\!0}+\alpha$ be a smooth irreducible flat $L$-connection, then from the previous consideration there exists $\mcal{S}_{\alpha,\veps} \subset \eus{H}^*\!(E,L)_\ell$ such that $p^*_\ell(\mcal{S}_{\alpha,\veps})$ is an open neighbourhood of $j_\ell([\nabla])$ in $\eus{M}^*\!(E,L)_\ell$. But from Lemma~\ref{Kuranishi smooth} we get $\mcal{S}_{\alpha,\veps} \subset \eus{H}^*\!(E,L)$, hence we have $p^*_\ell(\mcal{S}_{\alpha,\veps}) \subset j_\ell(\eus{M}^*\!(E,L))$, so we are done.}

\theorem{The moduli space $\eus{M}^*\!(E,L)$ of gauge equivalence classes of irreducible flat $L$-connections on $E$ has a structure of a topological space such that for each $[\nabla] \in \eus{M}^*\!(E,L)$ represented by $\nabla= \nabla_0 + \alpha \in \eus{H}^*\!(E,L)$ there exists an open neighbourhood $\mcal{U}_\alpha$ of $[\nabla]$ in $\eus{M}^*\!(E,L)$, an open neighbourhood $\mcal{O}_\alpha$ of $0$ in $\mcal{H}^1(E,\nabla)$ and a smooth mapping
\begin{align}
  \Phi \colon \mcal{O}_\alpha \rarr \mcal{H}^2(E,\nabla),
\end{align}
called the \emph{obstruction mapping}, satisfying $\Phi(0)=0$ and
\begin{align}
  \mcal{U}_\alpha \simeq \Phi^{-1}(0).
\end{align}
Thus $\mcal{U}_\alpha$ is homeomorphic to a closed subset in an open subset in a finite dimensional vector space.}

\proof{As the Kuranishi mapping $K_\nabla \colon \Omega^1_L(M,\End(E))_\ell \rarr \Omega^1_L(M,\End(E))_\ell$ is a local diffeomorphism at $0$ for $\ell > {1\over 2} \dim M +1$, there exist open neighbourhoods $\mcal{U},\mcal{V}$ of $0$ in $\Omega^1_L(M,\End(E))_\ell$ such that $K_{\nabla|\mcal{U}} \colon \mcal{U} \rarr \mcal{V}$ is a diffeomorphism of Hilbert manifolds.
We can take $\mcal{U}=\{\beta \in \Omega^1_L(M,\End(E))_\ell;\, \|\beta\|_\ell < \veps\}$ for $\veps > 0$ small enough, hence we have $\mcal{S}^\ell_\veps \subset \mcal{U}$. Furthermore, denote $F=(K_{\nabla|\mcal{U}})^{-1} \colon \mcal{V} \rarr \mcal{U}$. Because $\mcal{H}^1(E,\nabla)$ is a closed subspace of $\Omega^1_L(M,\End(E))_\ell$ and $\mcal{O} = \mcal{V} \cap \mcal{H}^1(E,\nabla)$ is an open subset of $\mcal{H}^1(E,\nabla)$, we obtain that $\mcal{O}$ is a Hilbert submanifold of $\Omega^1_L(M,\End(E))_\ell$. If we define the obstruction mapping $\Phi \colon \mcal{O} \rarr \mcal{H}^2(E,\nabla)$ by
\begin{align*}
\Phi(\gamma)= H_2\big([F(\gamma),F(\gamma)]\big),
\end{align*}
then $\Phi$ is a smooth mapping of Hilbert manifolds.

From the previous we have $K_\nabla(\mcal{S}^\ell_\veps) \subset \mcal{V} \cap \mcal{H}^1(E,\nabla) = \mcal{O}$. It remains to show that $K_\nabla(\mcal{S}^\ell_\veps)=\Phi^{-1}(0)$.
In case $\beta \in \mcal{S}^\ell_\veps$, then, using \eqref{decomposition}, we obtain $(\Phi \circ K_\nabla) (\beta)=H_2([\beta,\beta])=0$.
But on the other hand, if $\gamma \in \Phi^{-1}(0)$, then there exists a unique $\beta \in \mcal{U}$ satisfying $K_\nabla(\beta)=\gamma$. Hence we have $0=\Phi(\gamma)=(\Phi \circ K_\nabla)(\beta)=H_2([\beta,\beta])$. Since $\gamma \in \mcal{H}^1(E,\nabla)$, we get
\begin{align*}
  0 & = d^\nabla_1\!\gamma = d^\nabla_1\!(\beta +{\textstyle {1 \over 2}}( \delta^\nabla_1 \circ G_2)([\beta,\beta])), \\
  0 & = \delta^\nabla_0\! \gamma = \delta^\nabla_0\! \beta.
\end{align*}
Applying the Hodge decomposition \eqref{identity-flat} to the element ${1 \over 2} [\beta,\beta]$ and using the above equations, we obtain
\begin{align*}
  d^\nabla_1\! \beta + {\textstyle {1 \over 2}} [\beta,\beta] &= d^\nabla_1\! \beta + {\textstyle {1 \over 2}} (\delta^\nabla_2 \circ  d^\nabla_2 \circ G_2)([\beta,\beta])
   + {\textstyle {1 \over 2}}(d^\nabla_1 \circ  \delta^\nabla_1 \circ G_2)([\beta,\beta]) + {\textstyle {1\over 2}}H_2([\beta,\beta]) \\
  &={\textstyle {1 \over 2}}(\delta^\nabla_2 \circ  d^\nabla_2 \circ G_2)([\beta,\beta]) + {\textstyle {1\over 2}}H_2([\beta,\beta])={\textstyle {1 \over 2}}(\delta^\nabla_2 \circ  d^\nabla_2 \circ G_2)([\beta,\beta]).
\end{align*}
Denoting the left hand side of the equation above by $\psi$, we have
\begin{align*}
    \psi&=d^\nabla_1\! \beta + {\textstyle {1 \over 2}} [\beta,\beta] = {\textstyle {1 \over 2}}(\delta^\nabla_2 \circ  d^\nabla_2 \circ G_2)([\beta,\beta]) \\
    &={\textstyle {1 \over 2}}(\delta^\nabla_2 \circ G_3 \circ d^\nabla_2)([\beta,\beta])= {\textstyle {1 \over 2}}( \delta^\nabla_2 \circ G_3 )([d^\nabla_1\!\beta,\beta]- [\beta,d^\nabla_1\!\beta]) \\
    &= ( \delta^\nabla_2 \circ G_3 )([d^\nabla_1\!\beta,\beta])= (\delta^\nabla_2 \circ G_3)([\psi,\beta]),
\end{align*}
where we used that $[[\beta,\beta],\beta]=0$. Using the fact that there exists a positive constant $c$ such that
\begin{align*}
  \|(\delta^\nabla_2 \circ G_3)\varphi \|_\ell \leq c\, \|\varphi \|_{\ell-1},
\end{align*}
for all $\varphi \in \Omega^3_L(M,\End(E))_{\ell-1}$,
we make the following estimate
\begin{align*}
  \|\psi\|_{\ell-1}\leq \|\psi\|_\ell=\|(\delta^\nabla_2 \circ G_3)([\psi,\beta])\|_\ell \leq c\, \|[\psi,\beta]\|_{\ell-1} \leq c' \|\psi\|_{\ell-1} \|\beta\|_\ell < \veps c' \|\psi\|_{\ell-1},
\end{align*}
where $c'$ is another positive constant and the last inequality is provided that $\|\psi\|_{\ell-1} > 0$. If we take $\veps < {1 \over c'}$, then we have $\psi=0$.
Thus, together with $\delta^\nabla_0\!\beta =0$, we obtain that $\beta \in \mcal{S}^\ell_\veps$.
\smallskip

Further, because $j_\ell \colon \eus{M}^*\!(E,L) \rarr \eus{M}^*\!(E,L)_\ell$ is injective for $\ell > {1 \over 2} \dim M +2$, the mapping
\begin{align*}
  j_{k\ell|j_\ell(\eus{M}^*\!(E,L))} \colon j_\ell(\eus{M}^*\!(E,L)) \rarr j_k(\eus{M}^*\!(E,L))
\end{align*}
is bijective for $\ell \geq k >{1 \over 2} \dim M +2$ since $j_{k\ell} \circ j_\ell = j_k$. Moreover, form Lemma~\ref{open image} we know that $j_\ell$ has open image, therefore for each $\nabla_{\!0}+\alpha \in \eus{H}^*\!(E,L)$ there exists $\veps > 0$ satisfying that $p^*_\ell(\mcal{S}^\ell_{\alpha,\veps})$ is an open neighbourhood of $j_\ell([\nabla_{\!0}+\alpha])$ in $j_\ell(\eus{M}^*\!(E,L))$. Furthermore, from the previous we have that the following mapping
\begin{align*}
  p^*_\ell(\mcal{S}^\ell_{\alpha,\veps}) \xrarr{(p^\ell_{\alpha,\veps})^{-1}} \mcal{S}^\ell_{\alpha,\veps} \xrarr{\chi^\ell_\alpha} \mcal{S}^\ell_\veps \xrarr{K_\nabla} K_\nabla(\mcal{S}^\ell_\veps) \subset \mcal{O}^\ell_\veps = \mcal{V}^\ell_\veps \cap \mcal{H}^1(E,\nabla),
\end{align*}
where $\chi^\ell_\alpha \colon \mcal{S}^\ell_{\alpha,\veps} \rarr \mcal{S}^\ell_\veps$ is given through $\chi^\ell_\alpha(\nabla_{\!0}+\alpha+\beta)=\beta$, is a homeomorphism. Since $K_\nabla(\mcal{S}^\ell_\veps) \subset K_\nabla(\mcal{S}^k_\veps)$, for $\veps$ small enough we have the following commutative diagram
\begin{align*}
 \bfig
 \square<700,500>[p^*_\ell(\mcal{S}^\ell_{\alpha,\veps})`K_\nabla(\mcal{S}^\ell_\veps)`
 p^*_k(\mcal{S}^k_{\alpha,\veps})`K_\nabla(\mcal{S}^k_\veps);`j_{k\ell}`{\rm id}_{\mcal{H}^1(E,\nabla)}`]
 \efig
\end{align*}
in which ${\rm id}_{\mcal{H}^1(E,\nabla)}$ is a continuous mapping with respect to the norm $\|\cdot\|_k$ and $\|\cdot\|_\ell$ on $\mcal{H}^1(E,\nabla)$ because all norms on a finite dimensional vector space are equivalent. On the other hand, because we can find $\veps' \leq \veps$ such that $K_\nabla(\mcal{S}^k_{\veps'}) \subset K_\nabla(\mcal{S}^\ell_\veps)$, we obtain the following commutative diagram
\begin{align*}
 \bfig
 \square/>`<-`<-`>/<700,500>[p^*_\ell(\mcal{S}^\ell_{\alpha,\veps})`K_\nabla(\mcal{S}^\ell_\veps)`
 p^*_k(\mcal{S}^k_{\alpha,\veps'})`K_\nabla(\mcal{S}^k_{\veps'});`(j_{k\ell})^{-1}`{\rm id}_{\mcal{H}^1(E,\nabla)}`]
 \efig
\end{align*}
which gives that $j_{k\ell|j_\ell(\eus{M}^*\!(E,L))} \colon j_\ell(\eus{M}^*\!(E,L)) \rarr j_k(\eus{M}^*\!(E,L))$ is a homeomorphism.

Hence we have proved that $j_{k\ell|j_\ell(\eus{M}^*\!(E,L))} \colon j_\ell(\eus{M}^*\!(E,L)) \rarr j_k(\eus{M}^*\!(E,L))$ is a homeomorphism. Thus, $j_\ell$ gives a topology on $\eus{M}^*\!(E,L)$ which is independent on the Sobolev index $\ell$ for $\ell > {1 \over 2} \dim M +2$, and for each $\nabla= \nabla_{\!0} + \alpha$ there exists an open neighbourhood $\mcal{U}_\alpha = (j_\ell)^{-1}(p^*_\ell(\mcal{S}^\ell_{\alpha,\veps}))$ of $[\nabla]$ homeomorphic to $\Phi^{-1}(0)$.}

\theorem{The moduli space $\eus{M}^*\!(E,L)$ of irreducible flat $L$-connections on $E$ is a real (complex) manifold in a neighbourhood of $[\nabla] \in \eus{M}^*\!(E,L)$ if $H^2(E,\nabla)=0$, and its tangent space at $[\nabla]$ is naturally isomorphic to $H^1(E,\nabla)$.}

\proof{It follows immediately from the previous theorem.}

The moduli space $\eus{M}^*\!(E,L)$ is non-Hausdorff in general. The pairs of non-separable isomorphism classes of irreducible flat $L$-connections on $E$ can be detected by the following simple criterion.
\medskip

\lemma{If two distinct points $[\nabla]$ and $[\nabla']$ can not be separated in the moduli space $\eus{M}^*\!(E,L)$, then there exists a nontrivial homomorphism $\varphi \colon (E,\nabla) \rarr (E,\nabla')$ and $\psi \colon (E,\nabla') \rarr (E,\nabla)$ such that $\varphi \circ \psi=0$ and $\psi \circ \varphi=0$.}

\remark By a homomorphism $\varphi \colon (E,\nabla) \rarr (F,\nabla')$, where $\nabla$ and $\nabla'$ are $L$-connections on $E$ and $F$, respectively, we understand a homomorphism $\varphi \colon E \rarr F$ of vector bundles moreover satisfying the condition $\varphi \circ \nabla_{\!\xi} = \nabla'_{\!\xi} \circ \varphi$ for all $\xi \in \mfrak{X}(M)$.
\medskip

\proof{Because the moduli space $\eus{M}^*\!(E,L)$ is homeomorphic to the image of $j_\ell \colon \eus{M}^*\!(E,L) \rarr \eus{M}^*\!(E,L)_\ell$ and $\eus{H}^*\!(E,L)_\ell$ is a second-countable space for $\ell > {1\over 2} \dim M +2$, non-separability of $[\nabla]$ and $[\nabla']$ is equivalent to the existence of sequences $\nabla_{\!n}$ in $\eus{H}^*\!(E,L)_\ell$ and $\varphi_n$ in $\Gau(E)_{\ell+1}$ such that
\begin{align*}
  \nabla_{\!n} \rarr \nabla \quad \text{and} \quad \nabla'_{\!n} = \nabla_{\!n} \cdot \varphi_n \rarr \nabla'
\end{align*}
in $\eus{H}^*\!(E,L)_\ell$. Since $\nabla_{\!n} \cdot \varphi_n = \nabla'_{\!n}$, we obtain that $\varphi_n \in \ker \nabla''_{\!n}$, where $\nabla''_{\!n} \in \eus{H}(\End E,L)_\ell$ is the corresponding induced $L$-connection on $\End E$. Moreover, we have
\begin{align*}
  \nabla''_{\!n} \rarr \nabla'',
\end{align*}
where $\nabla'' \in \eus{H}(\End E,L)$ is the induce $L$-connection on $\End E$ given through $\nabla$ and $\nabla'$ on $E$.

Furthermore, because the real function $\nabla \mapsto \dim \ker \nabla$ is an upper semi-continuous function on $\eus{A}(\End E,L)_\ell$, we get
\begin{align*}
  \dim \ker \nabla'' \geq \limsup_{n \rarr +\infty} \dim \ker \nabla''_{\!n} \geq 1.
\end{align*}
Since $\nabla''$ is a smooth $L$-connection, there is a nonzero homomorphism $\varphi \colon (E,\nabla) \rarr (E,\nabla')$. By interchanging the role of $\nabla$ and $\nabla'$ in the argument above, we obtain a nonzero homomorphism $\psi \colon (E,\nabla') \rarr (E,\nabla)$.

Further, since $\nabla$ is an irreducible $L$-connection and $\psi \circ \varphi \in \ker \nabla^{\End E}$, we have $\psi \circ \varphi = c \cdot {\rm id}_E$ with $c \in \mathbb{K}$. As $\varphi \circ \psi = c \cdot {\rm id}_E$, it follows that $\varphi$ and $\psi$ are isomorphism if $c \neq 0$. Since $[\nabla] \neq [\nabla']$ by assumption, we get $c=0$.}

\section{Moduli spaces -- examples}

In this section we give some examples of moduli spaces of irreducible flat $L$-connections.
\medskip

\example (Holomorphic structures) Let $M$ be a complex manifold, i.e., a real manifold with a complex structure $\mcal{J}$, and consider the complex Lie algebroid $(L \rarr M,[\cdot\,,\cdot],a)$ defined in the following way. As a vector bundle $L=TM^{0,1}$, the anchor map $a \colon TM^{0,1} \rarr TM_\mathbb{C}$ is inclusion only and the Lie bracket is the Lie bracket of complexified vector fields.

Let $E \rarr M$ be a complex vector bundle, then an $L$-connection $\nabla$ on $E$ is a $\mathbb{C}$-linear mapping
\begin{align*}
  \nabla \colon \Gamma(M,E) \rarr \Gamma(M,T^*\!M^{0,1} \otimes E)
\end{align*}
which satisfies
\begin{align*}
  \nabla(fs)=d_Lf \otimes s + f \nabla s,
\end{align*}
where $s \in \Gamma(M,E)$ and $f \in C^\infty(M,\mathbb{C})$. Since we have $d_L=\bar{\partial} \colon C^\infty(M,\mathbb{C}) \rarr \Gamma(M,T^*\!M^{0,1})$, we obtain that a flat $L$-connection $\nabla$ on $E$ corresponds to a holomorphic structure $\bar{\partial}_E=\nabla$ on $E$.
\smallskip

Let $M$ be a connected compact Riemann surface $\Sigma_g$ of genus $g$. Then the deformation complex \eqref{deformation-complex} for a flat $L$-connection $\nabla$ is of the form
\begin{align*}
  0 \xrarr{\phantom{d^\nabla}} \Omega^0_L(M,\End E) \xrarr{d^\nabla} \Omega^1_L(M,\End E) \xrarr{\phantom{d^\nabla}} 0
\end{align*}
since $\rk L = 1$. Now, using the Atiyah--Singer index theorem, the index $\ind\mcal{E}(E,\nabla)$ of this complex is
\begin{align*}
  \ind\mcal{E}(E,\nabla)= \dim H^0(E,\nabla) - \dim H^1(E,\nabla) = -r^2(g-1),
\end{align*}
where $\rk E=r$. Moreover, if $\nabla$ is the irreducible $L$-connection, then $\dim H^0(E,\nabla)=1$ and we get $\dim H^1(E,\nabla)=1+r^2(g-1)$.
\medskip

\example (Higgs and co-Higgs bundles) Let $M$ be a complex manifold, i.e., a real manifold with a complex structure $\mcal{J}$, and consider the complex Lie algebroid $(L \rarr M,[\cdot\,,\cdot],a)$ defined in the following way. As a vector bundle $L=TM^{0,1} \oplus T^*\!M^{1,0}$, the anchor map $a \colon TM^{0,1} \oplus T^*\!M^{1,0} \rarr TM_\mathbb{C}$ is given via $a(X+\xi)=X$, where $X \in \Gamma(M,TM^{0,1})$ and $\xi \in \Gamma(M,T^*\!M^{1,0})$, while the Lie bracket is given through $[X+\xi,Y+\eta]_L=[X,Y]+\mcal{L}_X\eta-i_Y d\xi$.

Let $E \rarr M$ be a complex vector bundle, then an $L$-connection $\nabla$ on $E$ is a $\mathbb{C}$-linear mapping
\begin{align*}
  \nabla \colon \Gamma(M,E) \rarr \Gamma(M,L^*\otimes E) = \Gamma(M,T^*\!M^{0,1}\otimes E) \oplus \Gamma(M,TM^{1,0}\otimes E)
\end{align*}
which satisfies
\begin{align*}
  \nabla(fs)=d_Lf \otimes s + f \nabla s,
\end{align*}
where $s \in \Gamma(M,E)$ and $f \in C^\infty(M,\mathbb{C})$. For $d_L \colon C^\infty(M,\mathbb{C}) \rarr \Gamma(M,T^*\!M^{0,1}) \oplus \Gamma(M,TM^{1,0})$ we obtain that $d_Lf=\bar{\partial}f$. Further, we denote by $\bar{\partial}_E \colon \Gamma(M,E) \rarr \Gamma(M,T^*\!M^{0,1}\otimes E)$ the $\mathbb{C}$-linear mapping defined by $\bar{\partial}_E= {\rm pr}_1 \circ \nabla$ and by $\Phi \colon \Gamma(M,E) \rarr \Gamma(M,TM^{1,0} \otimes E)$ the $\mathbb{C}$-linear mapping defined through $\Phi = {\rm pr}_2 \circ \nabla$, where ${\rm pr}_1$ and ${\rm pr}_2$ are the projections onto $\Gamma(M,T^*\!M^{0,1}\otimes E)$ and $\Gamma(M,TM^{1,0}\otimes E)$, respectively. Moreover, we have $\bar{\partial}_E(fs)= \bar{\partial}f \otimes s + f \bar{\partial}_Es$ and $\Phi(fs)=f\Phi(s)$ for $f \in C^\infty(M,\mathbb{C})$ and $s \in \Gamma(M,E)$.

The flatness of an $L$-connection $\nabla$ on $E$ is equivalent to the following conditions. The mapping $\bar{\partial}_E \colon \Gamma(M,E) \rarr \Gamma(M,T^*\!M^{0,1} \otimes E)$ defines a holomorphic structure on $E$ and $\Phi \in \Gamma(M,TM^{1,0} \otimes \End E)$ is a holomorphic section, i.e., $\Phi \in H^0(M,TM \otimes \End E)$, satisfying $\Phi \wedge \Phi =0$, where $\Phi \wedge \Phi \in H^0(M,\Lambda^2 TM \otimes \End E)$. This is the so called \emph{co-Higgs bundle}, see \cite{Gualtieri2007}, \cite{Hitchin2010} and \cite{Rayan2010}.

In case there exists a holomorphic isomorphism between the holomorphic tangent bundle $TM$ and the holomorphic cotangent bundle $T^*\!M$ (e.g.~holomorphic symplectic manifolds, K3 surfaces, torus), then a flat $L$-connection $\nabla$ on $E$ can be viewed as a Higgs bundle $(\bar{\partial}_E,\Phi)$.
\smallskip

Let $M$ be a connected compact Riemann surface $\Sigma_g$ of genus $g$. Then the deformation complex \eqref{deformation-complex} for a flat $L$-connection $\nabla$ is of the form
\begin{align*}
  0 \xrarr{\phantom{d^\nabla}} \Omega^0_L(M,\End E) \xrarr{d^\nabla} \Omega^1_L(M,\End E) \xrarr{d^\nabla} \Omega^2_L(M,\End E) \xrarr{\phantom{d^\nabla}} 0
\end{align*}
since $\rk L = 2$. Now, using the Atiyah--Singer index theorem, the index $\ind\mcal{E}(E,\nabla)$ of this complex is
\begin{align*}
  \ind\mcal{E}(E,\nabla)= \dim H^0(E,\nabla) - \dim H^1(E,\nabla) +  \dim H^2(E,\nabla) = r^2(2g-2),
\end{align*}
where $\rk E=r$. Moreover, if $\nabla$ is the irreducible $L$-connection and $\nabla$ is a nonsingular point of $\eus{M}^*\!(E,L)$, then $\dim H^0(E,\nabla)=1$ and $\dim H^2(E,\nabla)=0$ and we get $\dim H^1(E,\nabla)=1-2r^2(g-1)$.

\providecommand{\href}[2]{#2}\begingroup\raggedright\endgroup

\end{document}